\documentclass[11pt]{amsart}

\usepackage{a4}
\usepackage{amsthm,amsmath,amssymb,amsfonts}
\usepackage[small]{caption}
\usepackage[normalem]{ulem}
\usepackage{color}
\usepackage{dsfont}
\definecolor{nb}{rgb}{.6,.176,1}
\definecolor{sienna}{rgb}{.92,.222,.176}
\definecolor{darkgreen}{rgb}{0,.5,0}

\DeclareMathOperator{\supp}{supp}

\DeclareMathOperator{\nod}{nod}
\DeclareMathOperator{\lf}{lf}
\newtheorem{definition}{Definition}[section]

\newtheorem{theorem}[definition]{Theorem}

\newtheorem{cor}[definition]{Corollary}
\newtheorem{lemma}[definition]{Lemma}
\newtheorem{proposition}[definition]{Proposition}
\newtheorem{example}[definition]{Example}
\newtheorem{remark}[definition]{Remark}

\newcommand{\smallcal}[1]{
	{\mathchoice{\scriptstyle}{\scriptstyle}{\scriptscriptstyle}{\scriptscriptstyle}\mathcal{#1}}}
\newcommand{\smallx}{\smallcal{X}}
\newcommand{\smally}{\smallcal{Y}}
\newcommand{\tree}{\smallx}
\newcommand{\treen}{\tree_n}
\newcommand{\treeR}[1][R]{\tree\restricted{#1}}
\newcommand{\treenR}[1][R]{(\treen)\restricted{#1}}

\newcommand{\E}{\mathcal{E}}
\newcommand{\Etran}{\E_\mathrm{trans}}

\newcommand{\TKal}{\mathcal{T}_{\scriptscriptstyle\mathrm{K}}}

\newcommand{\Prohorov}[1]{d^{#1}_\mathrm{Pr}}
\newcommand{\dPr}[1][]{\Prohorov{#1}}
\newcommand{\dPrT}[1][(T,r)]{\Prohorov{#1}}
\newcommand{\dPrTn}{\dPrT[(T_n,r_n)]}
\newcommand{\dH}[1][]{d_\mathrm{H}^{#1}}
\newcommand{\dGP}{d_\mathrm{GP}}
\newcommand{\dGHP}{d_\mathrm{GHP}}
\newcommand{\dsGHP}{d_\mathrm{sGHP}}
\newcommand{\da}[1][d]{{#1}^\#}
\newcommand{\daGP}{\dGP^\#}

\newcommand{\dasGHP}{\dsGHP^\#}

\newcommand{\comment}[1]{}

\newcommand{\mmspacefont}{\mathbb}
\newcommand{\Xmm}{\mmspacefont{X}}
\newcommand{\XHB}{\Xmm_\mathrm{HB}}
\newcommand{\Xc}{\Xmm_\mathrm{c}}
\newcommand{\Xfin}{\Xmm_\mathrm{fin}}
\newcommand{\Xprob}{\Xmm_1}
\newcommand{\K}{\mmspacefont{K}}
\newcommand{\XmH}{\mathfrak{X}}
\newcommand{\XmHc}{\XmH_c}

\newcommand{\XmHcsupp}{\XmHc^\mathrm{supp}}
\newcommand{\lmbsymb}{\mathfrak{m}}
\newcommand{\lmbfct}[2]{\lmbsymb_{#1}^{#2}}
\newcommand{\lmb}[1][R]{\lmbfct{\delta}{#1}}
\newcommand{\glmb}[1][\delta]{\lmbfct{#1}{}}
\newcommand{\deltah}{\hat{\delta}}

\newcommand{\half}{{\mathchoice{\textstyle}{}{}{}\tfrac12}}
\newcommand{\sm}{\smallskip}

\newcommand{\Wt}{\tilde{W}}

\newcommand{\To}{\Longrightarrow}

\newcommand{\tov}[1][]{\xrightarrow[#1]{vag}}

\newcommand{\ton}[1][n]{\tspace\displaystyle\mathop{\longrightarrow}_{\scriptscriptstyle#1\to\infty}\tspace}	
\newcommand{\Tno}[1][n]{\tspace\displaystyle\mathop{\Longrightarrow}_{\scriptscriptstyle#1\to\infty}\tspace}	

\newcommand{\eqlaw}{\overset{\CL}{=}}

\renewcommand{\P}{\mathbb{P}}

\newcommand{\be}{\begin{equation}}
\newcommand{\ee}{\end{equation}}

\newcommand{\eea}{\end{eqnarray}}
\newcommand{\bean}{\begin{eqnarray*}}
\newcommand{\eean}{\end{eqnarray*}}

\newif\ifpctex

\newcommand{\eps}{\varepsilon}

\newcommand{\N}{\mathbb{N}}
\newcommand{\R}{\mathbb{R}}
\newcommand{\Q}{\mathbb{Q}}

\newcommand{\CL}{\mathcal{L}}

\newcommand{\CM}{\mathcal{M}}
\newcommand{\M}{\CM}

\newcommand{\bew}[1]{\begin{equation*}\label{#1}}
\newcommand{\bea}[1]{\begin{eqnarray}\label{#1}}
\newcommand{\beL}[2]{\begin{lemma}[#2]\label{#1}}
\newcommand{\beD}[2]{\begin{definition}[#2]\label{#1}}
\newcommand{\beT}[2]{\begin{theorem}[#2]\label{#1}}
\newcommand{\beP}[2]{\begin{proposition}[#2]\label{#1}}
\newcommand{\beC}[2]{\begin{cor}[#2]\label{#1}}

\newcommand{\tno}{\ton}

\newcommand{\lemref}[1]{Lemma~\ref{Lem:#1}}
\newcommand{\defref}[1]{Definition~\ref{Def:#1}}
\newcommand{\exref}[1]{Example~\ref{Exp:#1}}
\newcommand{\propref}[1]{Proposition~\ref{P:#1}}
\newcommand{\thmref}[1]{Theorem~\ref{T:#1}}
\newcommand{\secref}[1]{Section~\ref{S:#1}}

\newcommand{\corref}[1]{Corollary~\ref{C:#1}}

\newcommand{\nbd}{\protect\nobreakdash-\hspace{0pt}}	

\newcommand{\vphi}{\varphi}

\newcommand{\restricted}[1]{{\mspace{-1mu}\upharpoonright}_{#1}}

\newcommand{\closure}[1]{\wtspace\overline{#1}\wtspace}
\newcommand{\Bcl}{\overline{B}}

\newcommand{\nlim}{\lim_{n\to\infty}}

\newcommand{\nliminf}{\liminf_{n\to\infty}}
\newcommand{\diffd}{\mathrm{d}}
\newcommand{\integralspace}{\/\mathchoice{\;}{\,}{\,}{}}		
\newcommand{\integral}[4]{\int_{#1}^{#2} \diffd#4 \integralspace#3 }	
\newcommand{\plainint}[2]{\int \diffd#2\integralspace #1 }	
\newcommand{\inta}[3]{\integral{#1}{}{#2}{#3}}				

\newcommand{\bintamx}[4]{\int_{#1} #2 \integralspace #3\(\diffd#4\)}
\newcommand{\wtspace}{\mathchoice{\,}{\,}{}{}}			
\newcommand{\tspace}{\mathchoice{\,}{}{}{}}			
\newcommand{\set}[2]{\{\wtspace #1 : #2 \wtspace\}}		
\newcommand{\bset}[2]{\bigl\{\, #1 : #2 \,\bigr\}}

\newcommand{\bigldelimiter}[1]{\mathchoice{\bigl#1}{\bigl#1}{{\textstyle#1}}{{\scriptstyle#1}}}
\newcommand{\bigrdelimiter}[1]{\mathchoice{\bigr#1}{\bigr#1}{{\textstyle#1}}{{\scriptstyle#1}}}
\renewcommand{\(}{\bigldelimiter(}
\renewcommand{\)}{\bigrdelimiter)}

\newenvironment{enproof}{\begin{proof}\begin{enumerate}}{\qedhere\end{enumerate}\end{proof}}

\newcommand{\Rtree}{$\R$\nbd tree}
\newcommand\define{\emph}
\newcommand{\distmat}[2][m]{\nu_{#1}(#2)}

\newcommand\munodT[1][T]{\mu^{\nod}_{#1}}
\newcommand\mudegT[1][T]{\mu^{\deg}_{#1}}
\newcommand\munodTn{\munodT[T_n]}
\newcommand\mudegTn{\mudegT[T_n]}
\newcommand{\lambdaT}[1][(T,r)]{\lambda_{#1}}
\newcommand{\lambdaTn}{\lambdaT[(T_n, r_n)]}

\newcommand{\und}{\quad\text{and\/}\quad}
\newcommand{\bcdot}{\boldsymbol{\cdot}}
\newcommand\TKalgeo{\TKal^{\scriptscriptstyle\mathrm{geom}(1/2)}}

\providecommand{\ifshorten}{\iffalse}
\ifshorten \renewcommand{\sm}{}\usepackage{a4wide} \fi

\numberwithin{equation}{section}

\title[Gap between Gromov-vague and Gromov-Hausdorff-vague topology]{The gap between Gromov-vague and
Gromov-Hausdorff-vague topology\footnote{Preprint of \emph{Stochastic Process.\ Appl.\ 126(9):2527--2553}}}

\author{Siva Athreya}
\address{Siva Athreya \\ Indian Statistical Institute
8th Mile Mysore Road  \\ Bangalore 560059, India.}
\email{athreya@isibang.ac.in}
\thanks{}

\author{Wolfgang L\"ohr}
\address{Wolfgang L\"ohr\\ Fakult\"at f\"ur Mathematik\\ Universit\"at Duisburg-Essen\\
Thea-Leymann-Str.\ 9\\ 45117 Essen, Germany}
\thanks{}
\email{wolfgang.loehr@uni-due.de}

\author{Anita Winter}
\address{Anita Winter\\ Fakult\"at f\"ur Mathematik\\ Universit\"at Duisburg-Essen\\
Thea-Leymann-Str.\ 9\\ 45117 Essen, Germany}
\thanks{}
\email{anita.winter@uni-due.de}


\begin{document}

\begin{abstract}
In \cite{AthreyaLoehrWinter} an invariance principle is stated for a class of strong Markov processes on tree-like metric measure spaces.
It is shown that if the underlying spaces converge Gromov vaguely, then the processes converge in the sense of finite dimensional distributions. Further,
 if the underlying spaces converge Gromov-Hausdorff vaguely, then the processes converge weakly in path space.
  In this paper we systematically introduce and study the Gromov-vague and the
  Gromov-Hausdorff-vague topology on the space of equivalence classes of metric boundedly
  finite measure spaces. The latter topology is closely related to the
  Gromov-Hausdorff-Prohorov metric which is defined
  on different equivalence classes of metric measure spaces.

  We explain the necessity of these two topologies via several examples, and close the gap between them.
  That is, we show that convergence in Gromov-vague topology implies
  convergence in Gromov-Hausdorff-vague topology if and only if the
  so-called lower mass-bound property is satisfied.
  Furthermore, we prove and disprove Polishness of several spaces of
  metric measure spaces in the topologies mentioned above (summarized in
  Figure~\ref{fig:Poltable}).

  As an application, we consider the Galton-Watson tree with critical offspring distribution of finite variance
  conditioned to not get extinct, and construct the so-called Kallenberg-Kesten tree as the weak limit in Gromov-Hausdorff-vague topology when the edge length are scaled down to go to zero.
\end{abstract}

\maketitle

\begin{quote}
{\tiny \tableofcontents }
\end{quote}

\section{Introduction}
\label{S:intro}
The paper introduces the {\em Gromov-vague} and the {\em Gromov-Hausdorff-vague} topology. These are two notions of convergence of (equivalence classes of) metric boundedly finite measure spaces.
These are ``localized'' versions of the
Gromov-weak topology and a topology closely related to the Gromov-Hausdorff-Prohorov topology on (equivalence classes of) metric finite measure spaces. \smallskip

\noindent
{\bf Gromov-weak convergence and sampling. } The Gromov-weak topology originates from the weak topology  in the space of probability measures on a fixed metric space.
It is an example of a topology which comes with a canonical family of measures and convergence determining test functions. That is,
given a complete, separable metric space, $(X,r)$, we denote by
${\mathcal M}_1(X)$ the space of all Borel probability measures on $X$
and by $\bar{{\mathcal C}}(X):=\bar{{\mathcal C}}_{\mathbb{R}}(X)$ the
space of bounded, continuous $\mathbb{R}$-valued functions. A sequence
of probability measures $(\mu_n)$ converges {\em weakly} to $\mu$ in ${\mathcal M}_1(X)$
(abbreviated $\mu_n\To\mu$), as $n\to\infty$, if and only if
$\int\mathrm{d}\mu_n\,f\to\int\mathrm{d}\mu\,f$ in $\mathbb{R}$, as
$n\to\infty$, for all $f\in\bar{{\mathcal C}}(X)$.

We wish to consider sequence of measures that live on different
spaces. In such a case an immediate analogue of bounded continuous
functions is not available. To still be in a position to imitate the
notion of weak convergence, we rely on the following useful fact: for
a sequence $(\mu_n)$ in ${\mathcal M}_1(X)$ and $\mu\in{\mathcal
  M}_1(X)$,
\begin{equation}
\label{d:001}
  \mu_n\Tno\mu\hspace{1cm}\mbox{ if and only if }\hspace{1cm}\mu_n^{\otimes\mathbb{N}}\Tno\mu^{\otimes\mathbb{N}}.
\end{equation}
Indeed, the {\em ``if''} direction follows by the fact that projections to a single coordinate are continuous.
The {\em ``only if''} direction follows as the set of bounded continuous functions $\varphi\colon
X^{\mathbb{N}}\to\R$ of the form $\varphi((x_n)_{n\in\mathbb{N}})=\prod_{i=1}^N\varphi_i(x_i)$ for some
$N\in\N$, $\varphi_i\colon X\to\R$, $i=1,\ldots,N$, separates points in $X^{\mathbb{N}}$ and is multiplicatively
closed (see, for example, \cite[Theorem~2.7]{Loehr} for an argument how to use \cite{LeCam1957} to conclude from
here that integration over such test functions is even convergence determining for measures on $\M_1(X)$).

Consider now the set of bounded continuous
functions $\varphi\colon X^{\mathbb{N}}\to\R$ of the following form
\begin{equation}
\label{d:002a}
   \varphi=\tilde{\varphi}\circ R^{(X,r)},
\end{equation}
where $R^{(X,r)}$ denotes the map that sends a vector $(x_n)_{n\in\mathbb{N}}\in X^\N$ to the matrix $(r(x_i,x_j))_{1\le i<j}\in\mathbb{R}^{\mathbb{N}\choose 2}_+$ of mutual distances, and
$\tilde{\varphi}\in\bar{{\mathcal C}}(\R_+^{{\mathbb{N}\choose 2}})$ depends on finitely many coordinates only. A \define{(complete, separable) metric measure space} $(X, r,\mu)$ consists of a complete, separable metric space $(X,r)$ and a Borel measure $\mu$ on $X$. Denote by $\Xprob$ the space of measure preserving isometry classes of metric spaces equipped with a Borel probability measure. Then for each representative $(X,r,\mu)$ of an isometry class $\smallx\in\Xprob$
the image measure $R^{(X,r)}_\ast\mu^{\otimes\mathbb{N}}=\mu^{\otimes\mathbb{N}}\circ (R^{(X,r)})^{-1}\in{\mathcal M}_1(\mathbb{R}_+^{{\mathbb{N}\choose 2}})$ is the same and is referred to as
the {\em distance matrix distribution $\nu^{\smallx}$} of $\smallx$. It turns out that if the distance matrix distributions of two metric measure spaces coincide, then the metric measure spaces fall into the same isometry class. This is known as Gromov's reconstruction theorem (compare \cite[Chapter~3$\tfrac{1}{2}$]{MR2000d:53065}), and suggests
to consider the {\em Gromov-weak topology}, which is the
topology induced by the set of functions of the form
\begin{equation}\label{d:002}
	\Phi(X,r,\mu) = \inta{X^\N}{\varphi}{\mu^{\otimes \N}} =
	\plainint{\tilde\vphi}{\nu^\smallx},
\end{equation}
where $\vphi$ is of the form \eqref{d:002a}. As this set is multiplicatively closed we can conclude once more
that it is also convergence determining for metric measure spaces on $\Xprob$.

The Gromov-weak topology on spaces of metric measure spaces, prescribed by test functions as in (\ref{d:002}), originates from the work of Gromov in the context of
metric geometry, where it is induced by so-called box metrics.  In
\cite{GrevenPfaffelhuberWinter2009} the Gromov-weak topology on
complete, separable metric measure spaces was reintroduced via
convergence of the functions of the form (\ref{d:002}), and metrized
by the so-called \emph{Gromov-Prohorov metric}.  Recently, in
\cite{Loehr}, it was shown that Gromov's box metric and the
Gromov-Prohorov metric are bi-Lipschitz equivalent.

Independently of Gromov's work, however, the idea of proving convergence of random
$0$-hyperbolic metric measure spaces (that means trees) via ``finite
dimensional distributions'', i.e., with the help of test functions of
the form (\ref{d:002}), has been used before in probability theory. As
the land mark we consider \cite[Theorem~23]{Aldous1993}, which states
Gromov-weak convergence of suitably rescaled Galton-Watson trees
towards the so-called Brownian continuum random tree (CRT), where the Galton-Watson trees are
associated with an offspring distribution of finite variance,
conditioned on a growing number of nodes and equipped with the uniform
distribution on its nodes.  Other results using test functions which
imitate sampling include
\cite[Theorem~4]{GrevenPfaffelhuberWinter2009}, where the so-called
$\Lambda$-coalescent tree is constructed as a Gromov-weak limit of
finite trees. Furthermore, distributions of finite samples from metric
measure spaces are used in hypothesis testing and for providing
confidence intervals in the field of topological data analysis (see,
for example, \cite{BlumbergGalMandellPancia,Carlson2014}).   \smallskip

\noindent
{\bf From Gromov-weak to Gromov-Hausdorff weak convergence. }
The following embedding result is known from
\cite[Lemma~5.8]{GrevenPfaffelhuberWinter2009}.  A sequence
$(\smallx_n)$ converges Gromov-weakly to $\smallx$ in $\Xprob$ if and
only if there is a complete, separable metric space $(E,d)$ such that
(representatives of) all $\smallx_n$ and $\smallx$ can be embedded
measure-preserving isometrically into $(E,d)$ in a way that the image
measures under the isometries converge weakly to the image limit
measure.  Using this embedding procedure, we can also define a
stronger topology: We say that a sequence $(\smallx_n)$ converges {\em
  Gromov-Hausdorff-weakly} to $\smallx$ in $\Xprob$ if and only if
there is a metric space $(E,d)$ such that we can do the above
embedding in a way that, additionally, the supports of the measures
converge in Hausdorff distance.

This topology is closely related to the one introduced under the name \emph{measured Hausdorff topology}
in \cite{Fukaya1987} in the context of studying the asymptotics of eigenvalues of the Laplacian on collapsing
Riemannian manifolds, and extended from compact to Heine-Borel measure spaces in \cite{KuwaeShioya2003}. The
difference to the Gromov-Hausdorff weak topology is that, instead of the supports, the whole spaces are required
to converge in Hausdorff metric topology. This leads to different equivalence classes, and the connection is
discussed extensively in \secref{GHvague}.
In probability theory, the measured (Gromov-)Hausdorff topology was reintroduced and further discussed in
\cite{EvansWinter2006,Miermont2009}, and recently extended in \cite{AbrahamDelmasHoscheit2013} to complete,
locally compact length spaces equipped with locally finite measures.
\smallskip

\noindent {\bf Verification of convergence. }
As for the Gromov-Hausdorff-weak topology no canonical family of convergence determining functions is available,
a key question is how to actually verify convergence in Gromov-Hausdorff-weak topology?
According to the definition, first an embedding of the whole sequence in to the same metric space must be provided.
For random forests there has been the tradition to encode them (if possible) as excursions on compact intervals,
and showing then convergence of the associated excursions in the uniform topology. As the map that sends an
excursion to a tree-like metric measure space is continuous with respect to the Gromov-(Hausdorff)-weak topology
(\cite[Proposition~2.9]{AbrahamDelmasHoscheit:exittimes}, \cite[Theorem~4.8]{Loehr}), convergence statements
obtained by re-scaling the associated excursions always imply convergence Gromov-Hausdorff-weakly.  This
approach has been successfully applied to branching forests with a particular offspring distribution (see, for
example, \cite{MR1954248,Duquesne2003,GrevenPopovicWinter2009}).
However, except for a few prototype models there is no obvious way to assign to a random graph model an
excursion coming from a Markov process.
In such a situation,  Gromov-Hausdorff convergence and Gromov-weak convergence are shown separately (for
example, \cite{HaasMiermont2012,CurienHaas2013,AddarioBroutinGoldschmidtMiermont}), or the scaling results are
stated either without the measure, using Gromov-Hausdorff convergence (for example,
\cite{LeGall2007,MarckertMiermont2011,HeWinkel}), or only in the  Gromov-weak topology (for example,
\cite{GrevenPfaffelhuberWinter2013}). \smallskip

\noindent {\bf Closing the gap. }
It is known that, if all considered metric measure spaces satisfy a (common) uniform volume doubling property,
then Gromov-weak and Gromov-Hausdorff-weak topology are the same (\cite[Corollary~27.27]{Villani2008}).
``Volume doubling'' is a standard property for Riemannian manifolds and regular, self-similar fractals.
It is quite restrictive for random spaces, such as random recursive fractals or, important for us, random
$\R$\nbd trees. In particular, Aldous's Brownian CRT almost surely does not have the doubling property, as can
be seen from the estimates in \cite[Theorem~1.3]{Croydon08} (see also \cite{DuquesneWang14} for stable L\'evy trees).

If the uniform volume doubling property fails, Gromov-Hausdorff-weak convergence is in general not implied by
Gromov-weak convergence.
The gap between Gromov-weak and Gromov-Hausdorff-weak topology, however,  sometimes matters a lot. \smallskip

\noindent {\bf Important example. }
We have recently considered in  \cite{AthreyaLoehrWinter} a class of strong Markov processes
on natural scale with values in $0$-hyperbolic compact metric spaces, which are uniquely determined by their
speed measures. We obtained  in \cite[Theorem~1]{AthreyaLoehrWinter} an invariance principle which states
convergence of such processes in path space provided the underlying  metric (speed-)measure spaces converge
Gromov-Hausdorff-weakly. If we only assume Gromov-weak convergence, the processes still converge in their finite
dimensional distributions, but without the additional convergence of the supports, convergence in paths space
fails. \smallskip

The main goal of the present paper is to close this gap between Gromov-Hausdorff-weak and Gromov-weak topology. We show that provided metric measure spaces converge Gromov-weakly, they also converge Gromov-Hausdorff-weakly if and only if the so-called (global) lower mass-bound property (Definition~\ref{Deef:001}) is satisfied.
This allows to verify Gromov-Hausdorff weak convergence via the following two steps (Theorem~\ref{T:mGH}):
\begin{enumerate}
\item Verify convergence of the test functions from (\ref{d:002}) together with
\item an extra ``tightness condition'' given by this lower mass-bound property.
\end{enumerate}
The same lower mass function also turns out to be useful for characterizing the metric measure spaces which are compact
and Heine-Borel, respectively, and for proving that the subspaces consisting of these metric measure spaces are Lusin
spaces but not Polish if equipped with the Gromov-weak topology.
The lower mass-bound property also appears in a compactness condition for the Gromov-Hausdorff-weak topology
(\corref{GHWcompact}).
Furthermore, we also extend the space of complete, separable metric probability measure spaces to complete,
separable, metric boundedly finite measure spaces and equip the latter with the so-called
Gromov-(Hausdorff)-vague topologies.

\subsection*{Outline}
The paper is organized as follows: In Section~\ref{S:Gromovvague} we recall the Gromov-weak topology on the space of metric finite measure spaces and then use it to define the Gromov-vague topology on the space of   metric boundedly finite measure spaces.
In Section~\ref{S:massbound} the global and local lower mass-bound properties  are defined and used to characterize compact metric (finite) measure spaces and Heine-Borel metric boundedly finite measure spaces.
In Section~\ref{S:embed} we characterize Gromov-vague convergence via isometric embeddings and deduce criteria for
Gromov-vague compactness and Gromov-vague tightness, as well as Polishness of the space of metric boundedly finite measure spaces
in Gromov-vague topology. Furthermore, we show that the subspaces of all compact and all Heine-Borel spaces,
respectively, are Lusin but not Polish.
In Section~\ref{S:GHvague} we introduce the stronger Gromov-Hausdorff-vague topology, and clarify its
relation to the measured Gromov-Hausdorff topology and the Gromov-Hausdorff-Prohorov metric. We also show that it is
a Polish topology on the space of Heine-Borel boundedly finite measure spaces. For the measured Gromov-Hausdorff
topology and the Gromov-Hausdorff-Prohorov metric, this means that restricting to spaces with measures of full support
yields again a Polish space.
In Section~\ref{S:relation} we prove our main convergence criterion for Gromov-Hausdorff-weak and -vague
convergence. Namely, given convergence in Gromov-weak or Gromov-vague topology, Gromov-Hausdorff-weak or
Gromov-Hausdorff-vague convergence is equivalent to the global or local lower mass-bound property, respectively.
In Section~\ref{S:excursions} we consider the construction of trees coded by continuous, transient excursions,
and show that the map which sends an excursion to the corresponding metric boundedly finite measure space is
continuous with respect to the Gromov-Hausdorff-vague topology.
Finally, as an example, we present the Gromov-Hausdorff-vague convergence in distribution of suitably re-scaled
finite-variance, critical Galton-Watson trees, which are conditioned on survival, to the so-called {\em continuum Kallenberg-Kesten tree}.

\section{The Gromov-vague topology}
\label{S:Gromovvague}
In this section we define the (pointed) Gromov-vague topology. We first introduce pointed metric boundedly finite
measure spaces, and the subspaces of interest. We recall the pointed Gromov-weak topology on pointed metric finite
measure spaces (Definition~\ref{Def:Gw}).  The pointed  Gromov-vague topology is then defined based on the Gromov-weak
topology via a ``localization procedure'' (Definition~\ref{Def:Gv}). We discuss the connection between both topologies
(Remark~\ref{Rem:002}), and present a perturbation result (Lemma~\ref{Lem:001}).

A \define{(pointed, complete, separable) metric measure space} $(X, r, \rho, \mu)$ consists of a complete, separable
metric space $(X,r)$, a distinguished point $\rho \in X$ called the \define{root}, and a Borel measure $\mu$ on $X$.
Since all our spaces are pointed, complete and separable, we usually drop these adjectives in the following when
referring to metric measure spaces.

\begin{definition}[equivalence of metric measure spaces]
Two metric measure spaces $(X,r,\rho,\mu)$ and $(X',r',\rho',\mu')$ are said to be \define{equivalent}
if and only if there is an isometry $\phi\colon \supp(\mu)\cup\{\rho\}\to\supp(\mu')\cup\{\rho'\}$ such that
$\phi(\rho)=\rho'$ and $\phi_\ast\mu=\mu'$, where as usual we denote by
\begin{equation}
\label{e:push}
   \phi_\ast\mu:=\mu\circ\phi^{-1}
\end{equation}
the \emph{push forward} of the measure $\mu$ under the measurable map $\phi$. We denote the equivalence of metric
measure spaces by $\cong$. Most of the time, however, we do not distinguish between a metric measure space and
its equivalence class.
\label{Def:000}
\end{definition}\sm

Recall that a \define{Heine-Borel space} is a metric space in which every bounded, closed set is compact. A
Heine-Borel space is obviously complete, separable and locally compact. We consider the following subclasses of
metric measure spaces.

\begin{definition}[$\Xmm$, $\XHB$, $\Xc$] $ $
\begin{enumerate}
	\item A metric measure space $(X, r, \rho, \mu)$ is called \define{boundedly finite}
		if the measure $\mu$ is finite on all bounded subsets of\/ $X$.
		Let\/ $\Xmm$ be the  set of (equivalence classes of) metric boundedly finite measure spaces.
	\item $\smallx \in \Xmm$ is called \define{Heine-Borel locally finite measure space} if the equivalence
		class contains a representative\/ $\smallx=(X,r,\rho,\mu)$ such that\/ $(X,r)$ is a Heine-Borel
		space. Denote the subspace of Heine-Borel spaces in $\Xmm$ by\/ $\XHB$.
	\item An equivalence
		class $\smallx \in \XHB$ is called \define{compact metric finite measure space} if it contains a representative\/ $\smallx=(X,r,\rho,\mu)$ such that $(X,r)$ is a compact space.
		Denote the subspace of compact spaces in $\XHB$ by\/ $\Xc$.
\end{enumerate}
\label{Def:001}
\end{definition}\sm

We illustrate this definition with an example which is useful for considering continuum limits of trees.
\begin{example}[locally compact geodesic spaces and $\R$-trees]
	Recall that a geodesic space is a metric space in which every two points are connected by an isometric
	path, i.e.\ a path with length equal to the distance between these points. A geodesic space is called
	\Rtree\ if there is, up to reparametrization, only one simple path between every pair of points.
	It is a classical fact that every complete, locally compact geodesic space is a Heine-Borel space.
	In particular, $\XHB$ contains the subclass of complete, locally compact\/ $\R$-trees with Radon measures.
\label{Exp:001}
\hfill$\qed$
\end{example}\sm

As every Heine-Borel space is locally compact, the local compactness assumption on the geodesic space is obviously essential.
The following remark discusses why the completeness assumption is important as well.
 \begin{remark}[non-complete spaces]
	We can allow also non-complete spaces as elements of\/ $\Xmm$ by identifying them with their respective
	completions. Note, however, that Radon measures on non-complete metric spaces are not boundedly finite in
	general.
	Consider for example the binary tree $T:=\{\rho\}\cup\bigcup_{n\in\mathbb{N}}\{0,1\}^n$ with edges connecting
	$w\in\{0,1\}^n$ with\/ $(w,0)\in\{0,1\}^{n+1}$ and $(w,1)\in\{0,1\}^{n+1}$, $n\in\mathbb{N}_0$, equipped with a
	metric determined by  $r(w,(w,0))=r(w,(w,1)):=c^{-n}$ if $w\in\{0,1\}^n$, for some $c\in[\tfrac{1}{2},1)$, and equipped with the
	length measure (see, Example~\ref{Exp:002} for a detailed definition).  The length measure is indeed a Radon
	measure as all compact subtrees are contained in a subtree spanned by finitely many vertices. On the other
	hand\/ $(T,r)$ is bounded, but the length measure is not finite.
	Thus non-complete, locally compact\/ $\R$-trees with a Radon measure are not elements of\/ $\Xmm$ in general.

	Moreover, non-complete, locally compact\/ $\R$-trees with a boundedly finite measure are not elements of\/
	$\XHB$ in general, as their completions do not need to be locally compact.
	Take for example $T:=(0,1] \times \{0\} \cup \bigcup_{n\in\N} \{\frac1n\} \times [0,1] \subseteq \R^2$, and let
	$r$ be the intrinsic length metric on $T$ (i.e., $r(x,y)$ is the Euclidean length of the shortest path within
	$T$ connecting $x$ and $y$). Then $(T,r)$ is a non-complete $\R$-tree, and it is easy to see that it is
	locally compact. Its completion $\overline{T}=T\cup\bigl\{(0,0)\bigr\}$, however, is not locally compact, because
	$(0,0)$ does not possess any compact neighborhood.
\label{Rem:001}
\hfill$\qed$
\end{remark}\sm

We next recall the definition of the (pointed) Gromov-weak topology on metric finite measure spaces (see
\cite{GrevenPfaffelhuberWinter2009} and \cite[Section~2.1]{LoehrVoisinWinter2013} for more details). As with the metric
measure spaces, we drop the adjective ``pointed'' in the following when referring to topologies on spaces of (pointed)
metric measure spaces.
\begin{definition}[(pointed) Gromov-weak topology]
	For $m\in\N$, the \define{$m$-point distance matrix distribution} of a metric finite measure space
	$\smallx=(X,r,\rho, \mu)$ is the finite measure on $\mathbb{R}_+^{\binom{m+1}2}$ defined by
	\begin{equation}
		\distmat{\smallx} := \int_{X^m}\mu^{\otimes
		m}(\mathrm{d}(x_1,\ldots,x_m))\,\delta_{(r(x_i,x_j))_{0\le i < j\le m}},
	\end{equation}
	where $x_0:=\rho$ and $\delta$ is the Dirac measure.
	A sequence $(\smallx_n)_{n\in\N}$ of metric finite measure spaces \define{converges} to a metric finite measure
	space $\smallx$ \define{Gromov-weakly} if all\/ $m$-point distance matrix distributions converge, i.e., if
	\begin{equation}\label{e:Gromovweak}
	  \distmat{\smallx_n} \Tno \distmat{\smallx},
	\end{equation}
	for all $m\in\mathbb{N}$, where we write $\Rightarrow$ for weak convergence of finite measures.
\label{Def:Gw}
\end{definition}\sm

Next we define the Gromov-vague topology on the space $\Xmm$ of metric boundedly finite measure spaces.
The construction is a straight-forward ``localization'' procedure, similar to the one used by Gromov for
Gromov-Hausdorff convergence of pointed locally compact spaces (compare \cite[Section~3B]{MR2000d:53065}).

Given a metric space $(X,r)$, we use the notations $B_r(x, R)$ and $\Bcl_r(x,R)$ for the open respectively closed
ball around $x\in X$ of radius $R\ge0$. If there is no risk of confusion, we sometimes drop the subscript $r$.
The restriction of a metric measure space $\smallx=(X, r, \rho, \mu)\in\Xmm$ to the closed ball $\Bcl(\rho, R)$
of radius $R\ge0$ around the root is denoted by
\begin{equation}
	\treeR := (X, r, \rho, \mu\restricted{\Bcl(\rho, R)})
		\cong \(\Bcl(\rho, R), r\restricted{\Bcl(\rho, R)^2}, \rho, \mu\restricted{\Bcl(\rho, R)}\).
\end{equation}

Generally (and informally), localization works as follows: given a topology on some class of spaces, the
localized form of convergence is defined for those spaces $\smallx$, where for all $R>0$, the
restriction $\smallx\restricted{R}$ falls into the original class.
Such spaces converge in the localized topology if, for almost all $R>0$, the restrictions converge.
If $d$ is a metric inducing the original topology, the localized convergence can therefore, for example, be
induced by the metric
\begin{equation}\label{e:da}
	\da\big(\smallx, \smally\big)
 :=
	\inta{\R_+}{e^{-R}\(1\land d(\smallx\restricted{R}, \smally\restricted{R})\)}R.
\end{equation}

We need the following lemma for our definition of Gromov-vague topology.
Denote the Gromov-Prohorov metric, which we define in Section~\ref{S:embed} below, by $\dGP$. For the moment, it
is enough to know that it induces the Gromov-weak topology by \cite[Theorem~5]{GrevenPfaffelhuberWinter2009}.

\begin{lemma}
	Let\/ $(\treen)_{n\in\N}$ be a sequence in\/ $\Xmm$ and\/ $\tree=(X, r, \rho, \mu)\in\Xmm$.
	The following are equivalent:
	\begin{enumerate}
		\item $\treenR\ton\treeR$ Gromov-weakly for all\/ $R>0$
			with\/ $\mu\(S_r(\rho, R)\) = 0$,  where  $S_r(\rho, R) = \Bcl_r(\rho, R)\setminus
			B_r(\rho,R)$ is the sphere of radius $R$ around $\rho$.
		\item $\treenR\ton\treeR$ Gromov-weakly for all but countably many\/ $R >0$.
		\item $\treenR \ton \treeR$ Gromov-weakly for Lebesgue-almost all\/ $R> 0$.
		\item There exists a sequence\/ $R_k\to\infty$ such that\/ $\treenR[R_k]\ton\treeR[R_k]$
			Gromov-weakly for all\/ $k\in\N$.
		\item $\daGP(\treen, \tree) \ton 0$.
	\end{enumerate}
\label{Lem:defGv}
 \end{lemma}

\begin{proof}
The implications {\em ``1.$\Rightarrow$2.$\Rightarrow$3.$\Rightarrow$4.''} are trivial. \sm

{\em ``4.$\Rightarrow$1.''} is a consequence of the
	Portmanteau theorem. Indeed, assume that $\treenR[R_k]\tno\treeR[R_k]$ Gromov-weakly along a sequence
	$R_k\to\infty$, and fix $R>0$. Choose $k\in\mathbb{N}$ large enough such that $R_k\ge R$.
	Then, for every $m\in\N$, $\distmat{\treenR[R_k]}\Tno \distmat{\treeR[R_k]}$.
	The first row of the $m$-point distance matrix $\nu_m$ contains, by definition, the distances to the
	root. Hence $\distmat{\treenR}$ is equal to the restriction of $\distmat{\treenR[R_k]}$ to the set of
	matrices with no entry in the first row exceeding $R$.
	The set of these matrices is closed, hence, by the Portmanteau theorem, the condition
	$\mu\(S_r(\rho, R)\) = 0$ implies the claimed convergence.

\sm\emph{``3.$\Leftrightarrow$5.''} follows directly from the fact that $\dGP$ induces the Gromov-weak topology,
	the definition of $\daGP$ in \eqref{e:da}, and the dominated convergence theorem.
\end{proof}\sm

We are now in a position to define the Gromov-vague topology.
\begin{definition}[(pointed) Gromov-vague topology]
	We say that a sequence $(\treen)_{n\in\N}$ in $\Xmm$ converges to $\tree\in\Xmm$ \define{Gromov-vaguely} if the
	equivalent conditions of \lemref{defGv} hold.
\label{Def:Gv}
\end{definition}\sm

Note that usually localized convergence is not strictly a generalization of the original one, because parts can
``vanish at infinity'' in the limit. For example, consider Gromov-Hausdorff convergence of (pointed)
compact metric spaces, and a sequence of two-point spaces, where the distance between the two points tends to
infinity. Such a sequence does not converge. In the localized Gromov-Hausdorff topology, however, it converges to
the compact space consisting of only one point. A similar phenomenon arises for the Gromov-vague topology.

\begin{remark}[Gromov-vague versus Gromov-weak]
	Consider the subspaces\/ $\Xfin$ and\/ $\Xprob$ of\/ $\Xmm$, consisting of spaces\/
	$\smallx = (X, r, \rho, \mu)$ where\/ $\mu$ is a finite measure, respectively a probability.
	Then on\/ $\Xprob$, the induced Gromov-vague topology coincides with the Gromov-weak topology.
	On\/ $\Xfin$, and even on\/ $\Xc$, however, this is
	not the case, because the total mass is not preserved in the Gromov-vague convergence.
	In fact, for\/ $\smallx=(X,r,\rho,\mu),\, \smallx_n=(X_n, r_n,\rho_n, \mu_n) \in \Xfin$,
	The Gromov-weak convergence $\smallx_n \to \smallx$ is equivalent to\/ $\smallx_n \to \smallx$
	Gromov-vaguely and\/ $\mu_n(X_n) \to \mu(X)$.
\label{Rem:002}
\hfill$\qed$
\end{remark}\sm

For a given metric space $(X,r)$, denote by $d_{\mathrm{Pr}}^{(X,r)}$ the Prohorov-metric on the space of all finite measures on $(X,{\mathcal B}(X))$,
i.e.,
\begin{equation}
\label{e:Pr}
\begin{aligned}
	&d_{\mathrm{Pr}}^{(X,r)}\big(\mu,\mu'\big)
  \\
 &:=
	 \inf\big\{\eps>0:\,\mu(A)\le\mu'(A^\eps)+\eps,\;\mu'(A) \le \mu(A^\eps)+\eps\;\;\forall A \text{ closed}\big\},
\end{aligned}
\end{equation}
where $A^\eps=\set{x}{d(x,A) \le \eps}$ is the closed $\eps$-neighborhood of $A$.
Recall that the Prohorov metric induces weak convergence.

We conclude this section with a simple stability property of Gromov-vague convergence under perturbations of the
measures in a localized Prohorov sense.
We will illustrate this later in Section~\ref{S:excursions} with Example~\ref{Exp:002}.

\begin{lemma}[perturbation of measures]
	Consider\/ $ \smallx = (X, r, \rho, \mu),\,\smallx_n = (X_n, r_n, \rho_n, \mu_n) \in \Xmm$, and
	another sequence of boundedly finite measure $\mu_n'$ on $X_n$, $n\in\N$. Assume that $\smallx_n \ton \smallx$
	Gromov-vaguely, and that there exists a sequence $R_k\to \infty$ such that for all $k\in\mathbb{N}$,
	\begin{equation}\label{Prmerg}
		\nlim d_{\mathrm{Pr}}^{(X_n,r_n)}\big(\mu_n\restricted{R_k}, \mu_n'\restricted{R_k}\big) =0.
	\end{equation}
	Then $\smallx_n^\prime := (X_n, r_n, \rho_n, \mu_n')$ converges Gromov-vaguely to $\smallx$.
\label{Lem:001}
\end{lemma}\sm

\begin{proof} Notice that for every fixed $k,n\in\mathbb{N}$,
\begin{equation}
\label{e:PrR}
    \lim_{R\downarrow R_k}d_{\mathrm{Pr}}^{(X_n,r_n)}\big(\mu_n\restricted{R},\mu_n\restricted{R_k}\big)=0.
\end{equation}
We may therefore assume w.l.o.g.\ that \eqref{Prmerg} and
$\treenR[R_k]\ton\treeR[R_k]$, Gromov-weakly, hold along the same sequence $(R_k)_{k\in\N}$. Thus for any fixed $k\in\mathbb{N}$,
\begin{equation}
\label{e:Gw001}
    (\smallx_n^\prime)\restricted{R_k} \ton \treeR[R_k],
\end{equation}
   Gromov-weakly, by Theorem~5 of \cite{GrevenPfaffelhuberWinter2009}. This, however, implies the claimed Gromov-vague convergence.
\end{proof}\sm

\section{The lower mass-bound property}
\label{S:massbound}
In this section we introduce the local and global lower mass-bound properties, and use them to
characterize compact spaces and Heine-Borel spaces, respectively.
These properties are formulated in terms of the following lower mass functions on the space of metric boundedly finite measure spaces.
For\/ $\delta,R>0$, we define $\lmb\colon \Xmm\to \R_+\cup\{\infty\}$ as
\begin{equation}
\label{lmb}
	\lmb\((X, r, \rho, \mu)\) := \inf\bset{\mu\(\Bcl_r(x,\delta)\)}{x\in B_r(\rho,R)\cap \supp(\mu)},
\end{equation}
with the convention that the infimum of the empty set is $\infty$ (which may happen if $\rho\not\in \supp(\mu)$).
Furthermore, set
\begin{equation}
\label{glmb}
	\glmb :=
    \lim_{R\to\infty} \lmb = \inf_{R > 0} \lmb.
\end{equation}

The following property plays an important r\^ole at several places in later arguments.
\begin{definition}[lower mass-bound property]
	A set\/ $\K\subseteq \Xmm$ of metric boundedly finite measure spaces
	satisfies the \define{local lower mass-bound property} if and only if
\begin{equation}
\label{e:lmb}
	   \inf_{\smallx \in\K} \lmb(\smallx)  > 0,
	\end{equation}
for all $R>\delta>0$. It satisfies the \define{global lower mass-bound property} if and only if
	\begin{equation}\label{e:glmb}
		\inf_{\smallx \in \K} \glmb(\smallx) > 0,
	\end{equation}
for all $\delta>0$.
We say that a single metric measure space $\smallx\in\Xmm$ satisfies the local/global mass-bound property if and
only if\/ $\mathbb{K}:=\{\smallx\}$ does.
\label{Deef:001}
\end{definition}\sm

Notice that in the definition of\/ $\lmb$, we could have replaced the closed ball by an open ball and/or the open
	ball by a closed ball without changing the conditions \eqref{e:lmb} and \eqref{e:glmb}.
	We made our choice such that\/ $\lmb$ is upper semi-continuous, which will be convenient
	in some proofs.

\begin{lemma}[upper semi-continuity]
	For every\/ $R,\delta>0$, the lower mass functions\/ $\lmb$ and\/ $\glmb$ are upper semi-continuous with
	respect to the Gromov-vague topology.
\label{Lem:lmbusc}
\end{lemma}\sm

\begin{proof} Fix $R,\delta>0$, and let $\smallx_n=(X_n,r_n,\rho_n,\mu_n) \to \smallx=(X,r,\rho,\mu)$ be a Gromov-vaguely
	converging sequence in $\Xmm$. Then we can choose $R'>R+\delta$ such that $\smallx_n' :=
	\smallx_n\restricted{R'}$ converges Gromov-weakly to $\smallx':=\smallx\restricted{R'}$.
	By Lemma~5.8 of \cite{GrevenPfaffelhuberWinter2009}, we can assume w.l.o.g.\ that $X$, $X_1$, $X_2$, ... are subspaces of some
	metric space $(E,d)$, and $\mu_n':= \mu_n\restricted{\Bcl(\rho,R')}$ converges weakly  to
	$\mu':=\mu\restricted{\Bcl(\rho, R')}$ on $(E,d)$.
	We can then find for every $x\in\supp(\mu)\cap B_r(\rho, R)$ a sequence $x_n\to x$ with $x_n\in\supp(\mu_n)$ for all $n\in\mathbb{N}$. Thus
\begin{equation}
\label{e:uppercon}
	\begin{aligned}
	\mu\(\Bcl(x, \delta)\)
 &=
   \inf_{\eps>0} \mu'\(\Bcl(x, \delta+\eps)\)
   \\
 &\ge
   \inf_{\eps>0} \nliminf \mu_n'\(\Bcl(x, \delta+\eps)\)
   \\
	 &\ge \nliminf \mu_n\(\Bcl(x_n,\delta)\)
 \\
  &\ge \lmb(\smallx_n),
	\end{aligned}
\end{equation}
	where we have applied the Portmanteau theorem in the second step, and used
in the last step that $x_n\in B(\rho_n, R)$ for large enough $n$.
	Hence $\lmb$ is upper semi-continuous. Therefore, $\glmb=\inf_{R>0}\lmb$ is also upper semi-continuous.
\end{proof}\sm

\begin{cor}[lower mass-bound property is preserved under closure]
	If\/ $\K\subseteq \Xmm$ satisfies the global or local lower mass-bound property, the same is true for its Gromov-vague
	closure\/ $\overline\K$.
\label{C:lmbcl}
\end{cor}\sm

\begin{lemma}[characterization of compact mm-spaces]
	Let\/ $\smallx\in\Xmm$. Then\/ $\smallx$ is a \emph{compact} metric finite measure space if and only if
	it has finite total mass, and satisfies the global lower mass-bound property.
\label{Lem:cchar}
\end{lemma}\sm

\begin{proof} {\em ``$\Rightarrow$'' } Assume that $\smallx=(X,r,\rho, \mu)$ is compact. Then $X$ is bounded,
	and hence $\mu$ is a finite measure.
	For every $\delta>0$, the function $x\mapsto \mu\(B(x, \delta)\)$ is lower semi-continuous. Therefore,
	it attains its minimum on the compact set $\supp(\mu)$, and thus the global lower mass-bound property holds.\sm

{\em ``$\Leftarrow$'' }
 Assume that $\mu$ is finite, and that the  global lower mass-bound property holds.
 Then for all $\delta>0$, we can
	cover $\supp(\mu)$ with finitely many balls of radius $2\delta$. To see this, notice that we can choose an at most countable covering $\{B(x,2\delta);\,x\in S\subseteq X\}$ of $\supp(\mu)$
with the property that the points in $S$ have mutual distances at least $2\delta$. As  $\{B(x,\delta);\,x\in
S\subseteq X\}$ then consists of pairwise disjoint sets, each carrying $\mu$-mass
at least $\glmb(\smallx)$,
the total mass of $\mu$ is at least $\glmb(\smallx)\cdot \# S$. As $\mu$ is a finite measure,
$\{B(x,2\delta);\,x\in S\subseteq X\}$ must be a finite set.
Since  $\supp(\mu)$ is complete, this means that $\supp(\mu)$ is actually compact.
\end{proof}\sm

\begin{lemma}[characterization of Heine-Borel mm-spaces]
	Let\/ $\smallx\in\Xmm$. Then $\smallx$ is a \emph{Heine-Borel} locally finite measure space if and only
	if it satisfies the local lower mass-bound property.
\label{Lem:HBchar}
\end{lemma}\sm

\begin{proof}
	Given $R>0$, $\smallx$ satisfies $\lmb(\smallx)>0$ for every $\delta>0$ if and only if $\smallx\restricted{R}$
	satisfies the global lower mass-bound property. Hence by \lemref{cchar}, $\smallx\restricted{R}$ satisfies the global lower mass-bound property if and only if
	$\smallx\restricted{R}$ is compact. Obviously, $\smallx\restricted{R}$ is compact for all $R>0$ if and only if
$\smallx$ is Heine-Borel.
\end{proof}\sm

\begin{cor}[$\XHB$ and $\Xc$ are measurable]
	Both\/ $\XHB$ and\/ $\Xc$ are measurable subsets of\/ $\Xmm$ with respect to Borel $\sigma$-field generated by the Gromov-vague topology.
\label{C:measurable}
\end{cor}\sm

\begin{proof} Notice that
\begin{equation}
\label{e:ARdeltaa}
   \XHB=\bigcap\nolimits_{R\in\N} \bigcap\nolimits_{\delta>0} \bigcup\nolimits_{a>0}\big\{\smallx\in\Xmm:\,\lmb(\smallx) \ge a\big\},
\end{equation}
by \lemref{HBchar}. Since the lower mass functions are upper semi-continuous by \lemref{lmbusc},
$A_{R,\delta,a}:=\set{\smallx\in\Xmm}{\lmb(\smallx) \ge a}$ is closed for all $\delta,R>0$.
Hence $\XHB$ is measurable. The measurability of $\Xc$ follows analogously by noticing that
\begin{equation}
\label{e:Adeltaa}
   \Xc=\bigcap\nolimits_{\delta>0} \bigcup\nolimits_{a>0}\set{\smallx\in\Xmm}{\glmb(\smallx) \ge a,\; \mu(X) \le
   a^{-1}},
\end{equation}
by \lemref{cchar}.
\end{proof}\sm

\section{Embeddings, compactness and Polishness}
\label{S:embed}
Recall that weak convergence of finite measures on a complete, separable metric space is
induced by the complete Prohorov metric (see, (\ref{e:Pr})). In the same spirit,
the Gromov-weak topology is induced by the
complete \emph{Gromov-Prohorov metric}, which is defined for two metric finite measure spaces $\smallx = (X,r,\mu)$ and $\smallx' = (X',r',\mu')$ by
\begin{equation}
\label{e:dGP}
	\dGP(\smallx, \smallx') := \inf_{d}d_{\mathrm{Pr}}^{(X\sqcup X',d)}\big(\mu,\mu'\big),
\end{equation}
where the infimum is taken over all metrics $d$ on $X\sqcup X'$ that extend both $r$ and $r'$, and $\sqcup$
denotes the disjoint union (see \cite[Theorem~5]{GrevenPfaffelhuberWinter2009}).

The fact that $\dGP$ induces the Gromov-weak topology immediately implies the following embedding result: for every Gromov-weakly convergent
sequence, $((X_n, r_n, \mu_n))_{n\in\mathbb{N}}$, there exists a common complete, separable metric space $(E,d)$ in which all
$(X_n, r_n)$ can be isometrically embedded such that (the push-forwards of) the measures $\mu_n$ converge weakly
to a measure $\mu$ on $(E,d)$ (compare \cite[Lemma~5.8]{GrevenPfaffelhuberWinter2009}).

In this section we show that an analogous statement (\propref{embed}) is true for the
Gromov-vague topology and, if the sequence satisfies the local lower mass-bound, $(E,d)$ can be chosen as Heine-Borel
space. We will apply this to characterize compact sets in $\mathbb{X}$ (Corollary~\ref{C:compact}), and to show that $\mathbb{X}$ is Polish (Proposition~\ref{P:Polish}), while
$\Xc$ and $\XHB$ are Lusin spaces but not Polish (Corollary~\ref{C:Lusin}).
We also prove a tightness criterion for probability measures on $\mathbb{X}$ (Corollary~\ref{C:tight}). \sm

We start with the embedding result.
\begin{proposition}[characterization via isometric embeddings]
	For each $n\in\N\cup\{\infty\}$, let\/ $\treen=(X_n, r_n, \rho_n, \mu_n)\in\Xmm$. Then $(\treen)_{n\in\N}$
	converges to $\tree_\infty$ Gromov-vaguely if and only if there exists a pointed complete, separable
	metric space $(E, d, \rho)$ and isometries $\varphi_n \colon \supp(\mu_n) \to E$ such that\/ $\varphi_n(\rho_n)
	= \rho$ for $n\in\N\cup\{\infty\}$, and
	\begin{equation}\label{e:embed}
	 	\((\varphi_n)_\ast\mu_n\)\restricted{\Bcl_d(R, \rho)} \,\Tno\,
		\((\varphi_\infty)_\ast\mu_\infty\)\restricted{\Bcl_d(R, \rho)},
	\end{equation}
for all but countably many $R\ge0$.
	Furthermore, if\/ $\set{\smallx_n}{n\in\N}$ satisfies the local lower mass-bound property,
	then $\smallx_\infty\in\XHB$ and\/ $E$ can be chosen as Heine-Borel space.
	In this case, \eqref{e:embed} is equivalent to
	\begin{equation}\label{e:embedv}
	 	(\varphi_n)_\ast\mu_n \,\tov\, (\varphi_\infty)_\ast\mu_\infty,
	\end{equation}
	where $\tov$ denotes vague convergence of Radon measures on $E$.
\label{P:embed}
\end{proposition}\sm

Before we give the proof, we illustrate with an example that the local lower mass-bound property cannot be
dropped without replacement in the second part of the proposition, even if the limit is assumed to be compact.
\begin{example}[$E$ is not Heine-Borel without lower mass-bound]
	Consider\/ $\smallx_n=\([0,1]^n, r_n, 0, \tfrac{n-1}{n}\delta_{0^n}+\tfrac{1}{n}\lambda_n\)$, where $r_n$ is the
	Euclidean metric, $\delta_{0^n}$ is the Dirac measure in\/ $0^n=(0,\ldots,0)\in [0,1]^n$, and\/
	$\lambda_n$ is the $n$\nbd dimensional Lebesgue measure.
	Then\/ $\smallx_n$ is compact and obviously converges Gromov-vaguely (and Gromov-weakly) to the
	compact probability space consisting of only one point, but the embedding space\/ $(E,d)$ cannot be
	chosen as Heine-Borel space.
\label{ex:conv}
\hfill$\qed$
\end{example}\sm

\begin{proof}[Proof of Proposition~\ref{P:embed}]
	It is easy to see that \eqref{e:embed} implies the Gromov-vague convergence, and that if $E$ is a
	Heine-Borel space, \eqref{e:embed} is equivalent to \eqref{e:embedv}.

	Conversely, assume that $\smallx_n\tno\smallx$  Gromov-vaguely, and abbreviate $X:=X_\infty$, $r:=r_\infty$ and
	$\rho:=\rho_\infty$.
	Let $(R_k)_{k\in\N}$ be an increasing sequence of radii with $\lim_{k\to\infty} R_k=\infty$ and
	$\mu\(\Bcl_r(\rho,R_k)\setminus B_r(\rho, R_k)\)=0$. Using that the Gromov-Prohorov metric metrizes the
	Gromov-weak topology by \cite[Theorem~5]{GrevenPfaffelhuberWinter2009}, we can construct for $n,k\in\N$
	a metric $d_{n,k}$ on $X_n \sqcup X$ extending both $r_n$ and $r$ such that for all $l\in\{1,\ldots,k\}$,
	\begin{equation}
\label{e:x00}
		\nlim d_{\mathrm{Pr}}^{(X_n \sqcup X,d_{n,k})}\big(\mu_n\restricted{R_l}, \mu\restricted{R_l}\big) = 0,
	\end{equation}
	where we use the abbreviation
	\begin{equation}
   \label{e:x01}
		\mu\restricted{R} := \mu\restricted{\Bcl_d(\rho, R)}.
	\end{equation}
	It is easy to check that we can do it such that $\rho_n$ and $\rho$ are identified.
	Using Cantor's diagonal argument, we can find a subsequence $(k_n)$ such that $d_n:=d_{n,k_n}$ satisfies
	$\nlim d_{\mathrm{Pr}}^{(X_n \sqcup X,d_{n})}(\mu_n\restricted{R_k}, \mu\restricted{R_k}) = 0$, for every $k\in\N$.
	Let $E:=\bigsqcup_{n\in\N \cup \{\infty\}} X_n$, $d$ the largest metric on $E$ which extends all $d_n$,
	and $\varphi_n\colon X_n \to E$ the canonical injection. Then it is easy to check that $(E,d)$ is a
	complete, separable metric space and \eqref{e:embed} is satisfied. \sm

	Now assume that $\set{\smallx_n}{n\in\N}$ satisfies the  local lower mass-bound property. Then it is also
	satisfied for $\bset{\smallx_n}{n\in\N \cup \{\infty\}}$ by \corref{lmbcl}.
	Due to \lemref{HBchar}, we may assume that $X_n$ and $X$ are Heine-Borel spaces. We
	have to show that $E$ is a Heine-Borel space as well.
	To this end, we show that every bounded sequence $(x_i)_{i\in\N}$ in $E$ has an accumulation point.
	If infinitely many $x_i$ are in a single $X_n$, $n\in\N\cup\{\infty\}$, this follows from the
	Heine-Borel property of $X_n$.
	Therefore, we can assume w.l.o.g.\ that $x_n \in X_n \cap B_d(\rho, R_k)$ for all $n$ and some $k$.
	By \eqref{e:lmb} together with $(\mu_n)\restricted{R_k} \Rightarrow \mu\restricted{R_k}$
	on $E$, we obtain $d(x_n, X) \to 0$. Hence there is $y_n \in X$ with $d(x_n, y_n)\to 0$ and, by the
	Heine-Borel property of $X$, $(y_n)_{n\in\N}$ has an accumulation point, which is also an accumulation
	point of $(x_n)_{n\in\N}$.
\end{proof}\sm

From here we can easily characterize the relatively compact sets.
\begin{cor}[Gromov-vague compactness]
	For a set\/ $\K\subseteq \Xmm$ the following are equivalent:
	\begin{enumerate}
	\item\label{it:Gv} $\K$ is relatively compact in $\Xmm$ equipped with the Gromov-vague topology.
	\item\label{it:Gw} For all\/ $R>0$, the set of restrictions $\K\restricted{R} :=
		\bset{\treeR}{\tree\in \K}$ is relatively compact in the Gromov-weak topology.
	\item\label{it:Gws} $\K\restricted{R_k}$ is relatively compact in the Gromov-weak topology for a
		sequence $R_k\to\infty$.
	\end{enumerate}
	Furthermore, a set\/ $\K\subseteq \XHB$ which satisfies the local lower mass-bound property is relatively compact in $\XHB$
    equipped with Gromov-vague topology if and only if the total masses of
	large balls are uniformly bounded, i.e., for all $R>0$,
	\begin{equation}\label{e:umb}
		\sup_{(X,r,\rho,\mu)\in\K} \mu\(B_r(\rho, R)\) < \infty.
	\end{equation}
\label{C:compact}
\end{cor}\sm

\begin{remark}[Gromov-weak compactness]
	Criteria for relative compactness in the Gromov-weak topology are given in Theorem~2 and Proposition~7.1 of\/
		\cite{GrevenPfaffelhuberWinter2009}.
\label{Rem:003}
\hfill$\qed$
\end{remark}\sm

\begin{remark}[convergence without the lower-mass bound property]
	As we have seen in Example~\ref{ex:conv}, a Gromov-vaguely convergent sequence in $\XHB$ does not have to satisfy the
	local lower mass-bound property. Hence the local lower mass-bound property is not necessary for
	relative compactness in $\XHB$.
\label{Rem:011}
\hfill$\qed$
\end{remark}\sm

\begin{proof}[Proof of Corollary~\ref{C:compact}]
{\em ``\ref{it:Gv}$\Rightarrow$\ref{it:Gw}'' } Assume that $\K$ is relatively compact.
	Let $R>0$, and consider a sequence $\(\smallx_n=(X_n, r_n, \rho_n, \mu_n)\)_{n\in\N}$
	in $\K$.  Then it possesses a Gromov-vague limit point $\smallx \in \Xmm$, and, by
	passing to a subsequence, we may assume w.l.o.g.\ that $\smallx_n\tno\smallx$ Gromov-vaguely.
	By \propref{embed}, we can assume that $X_n\subseteq E$ for some separable metric space $E$, and
	$\mu_n':=\mu_n\restricted{B(\rho, R')} \Tno \mu\restricted{B(\rho, R')} =: \mu'$ for some $R'> R$.
	By the Prohorov theorem, the sequence $(\mu_n')_{n\in\N}$ is tight. Thus
	$\(\mu_n\restricted{\Bcl(\rho, R)}\)_{n\in\N}$ is also tight.  We can conclude once more with the Prohorov theorem
	that $\(\mu_n\restricted{\Bcl(\rho, R)}\)_{n\in\N}$  is relatively weakly compact.
	Consequently, $(\smallx_n\restricted{R})_{n\in\N}$ has a Gromov-weak limit point.\sm

{\em ``\ref{it:Gw}$\Rightarrow$\ref{it:Gws}'' } is obvious. \sm

{\em ``\ref{it:Gws}$\Rightarrow$\ref{it:Gv}'' }
	Let $(\smallx_n)_{n\in\N}$ be a sequence in $\K$. By passing to a subsequence, we may assume that
	$\smallx_n\restricted{R_k}$ converges Gromov-weakly to some metric finite measure space $\smallx^{(k)}$
	for all $k$. Now it is easy to check that $\smallx^{(i)} = \smallx^{(k)}\restricted{R_i}$ whenever $R_i
	\le R_k$ and that we can therefore construct $\smallx\in\Xmm$ with
	$\smallx\restricted{R_k}=\smallx^{(k)}$ for every $k\in\N$. By definition, $\smallx_n \to \smallx$ in
	Gromov-vague topology.\sm

	Now assume that $\K\subseteq\XHB$ satisfies the local lower mass-bound property and \eqref{e:umb}. Fix $R>0$. Then for
	every $\eps>0$ we can find $N=N(\varepsilon,\K)\in\N$ such that for every $\smallx=(X,r,\rho, \mu)\in\K$, we can cover
	$\Bcl_r(\rho, R)$ by $N$ balls of radius $\eps$. Hence $\K\restricted{R}$ is relatively compact in
	Gromov-weak topology by Proposition~7.1 of \cite{GrevenPfaffelhuberWinter2009}. Therefore, $\K$ is
	relatively compact in $\Xmm$ with Gromov-vague topology. As also $\closure{\K}$ satisfies the local lower
	mass-bound property by \corref{lmbcl}, $\closure{\K}\subseteq \XHB$ by \lemref{HBchar}.
\end{proof}\sm

Having a characterization of compactness at hand, we can also characterize tightness of probability measures on $\Xmm$.
Denote by $\Xfin$ the subspace of metric finite measure spaces.
\begin{cor}[tightness of measures on $\Xmm$]
	Let\/ $\Gamma$ be a family of probability measures on $\Xmm$, and consider for each $R>0$ the restriction map
	$\psi_R \colon \Xmm \to \Xfin$ given by $\smallx \mapsto \smallx\restricted{R}$.
	Then the following are equivalent:
	\begin{enumerate}
	\item $\Gamma$ is Gromov-vaguely tight.
	\item The set $(\psi_R)_\ast(\Gamma)$ is Gromov-weakly tight for all $R>0$.
	\item For all $R, \eps >0$, there is a $\delta>0$ such that
	\begin{gather}\label{eq:modmass}
	    \sup_{\mathbb{P} \in \Gamma}
	      \mathbb{P}\big\{(X, r, \rho, \mu)\in \Xmm:\,\mu\big(B_r(\rho, R)\big) > \tfrac1\delta\big\} \le \eps,\\
	    \sup_{\P\in\Gamma}
	      \mathbb{P}\big\{(X,r,\rho,\mu)\in \Xmm:\,\mu\big\{x\in B_r(\rho, R):\,\mu(B_r(x,\eps))\le\delta\big\} \ge \eps\big\}
	      \le \eps. \nonumber
	\end{gather}
\comment{	\begin{equation}
    \label{eq:modmass}
	\begin{aligned}
	\sup_{\mathbb{P} \in \Gamma}\Big(
		&\mathbb{P}\big\{(X, r, \rho, \mu)\in \Xmm:\,\mu\big(B_r(\rho, R)\big) > \tfrac1\delta\big\}
  \\
		&{} + \mathbb{P}\big\{(X,r,\rho,\mu)\in \Xmm:\,\mu\big\{x\in B_r(\rho, R):\,\mu(B_r(x,\eps))\le\delta\big\} \ge \eps\big\}
			\Big)\le \eps.
	\end{aligned}
	\end{equation}
}
\comment{	\begin{equation}
		\sup_{P\in \Gamma}\, \bintamx\Xmm{ \mu\(\bset{x\in B_r(\rho, R)}{\mu(B_r(x,\eps)) \le \delta}\)
		}P{(X, r, \rho, \mu)} \le \eps.
	\end{equation}}
	\end{enumerate}
\label{C:tight}
\end{cor}\sm

\begin{remark}[Gromov-weak tightness] A characterization of Gromov-weak tightness of probability measures of metric finite measures spaces is given in  Theorem~3 in
	\cite{GrevenPfaffelhuberWinter2009} (compare also \cite[Remark~7.2(ii)]{GrevenPfaffelhuberWinter2009}).
\label{Rem:006}
\hfill$\qed$
\end{remark}\sm

\begin{proof}[Proof of Corollary~\ref{C:tight}] {\em ``only if'' } Assume that the family $\Gamma$ is
	Gromov-vaguely tight. Then we can find for all $\eps>0$ a compact set
	$\K_\eps$ with $\mathbb{P}(\K_\eps)\ge 1-\eps$ for all $\P\in\Gamma$. In particular, by
	\corref{compact}, the sets  $\K_\eps\restricted{R}$  are Gromov-weakly relatively compact
	 for all $R>0$. Because $(\psi_R)_*(\mathbb{P})(\K_\eps\restricted{R}) \ge \mathbb{P}(\K_\eps) \ge
	 1-\eps$ for all $\mathbb{P}\in\Gamma$, the set $(\psi_R)_*(\Gamma)$ is Gromov-weakly tight.
	\sm

	\emph{``if'' } Conversely if, for all $\eps,R>0$, $\K_{\eps, R}$ are Gromov-weakly compact sets satisfying
	$(\psi_R)_*(\mathbb{P})(K_{\eps,R}) \ge 1-\eps$,
	then for all $\eps>0$, $\K_\eps := \bset{\smallx\in \Xmm}{\smallx\restricted{n} \in\K_{2^{-n}\eps, n}\;\forall n\in\N}$ is
	a Gromov-vaguely relatively compact set which satisfies $\mathbb{P}(K_\eps)\ge 1-\eps$ for all $\mathbb{P}\in\Gamma$.\sm
	\sm

	The equivalence of \eqref{eq:modmass} now follows from Theorem~3 in
	\cite{GrevenPfaffelhuberWinter2009}.
\end{proof}\sm

Constructing a complete metric on $\Xmm$ that metrizes the Gromov-vague topology is now standard.

\begin{proposition}[$\Xmm$ is Polish]
	The space $\Xmm$ of metric boundedly finite measure spaces equipped with the Gromov-vague topology
	is a Polish space.
\label{P:Polish}
\end{proposition}\sm

\begin{proof}
	One possible choice of a complete metric is
	\begin{equation}
		\daGP\big(\smallx, \smally\big)
 :=
		\inta{\R_+}{e^{-R}\(1\land \dGP(\smallx\restricted{R}, \smally\restricted{R})\)}R.
	\end{equation}
	Indeed, that $\daGP$ induces the Gromov-vague topology is shown in Lemma~\ref{Lem:defGv}, and
	separability is obvious.
	To see completeness, consider a Cauchy sequence $(\smallx_n)_{n\in\N}$ in $\Xmm$.
	Then $\smallx_n\restricted{R}$ is a Cauchy sequence with respect to $\dGP$ for Lebesgue-almost all $R>0$.
	By completeness of $\dGP$, $\set{\smallx_n\restricted{R}}{n\in\N}$ is relatively compact in the
	Gromov-weak topology for these $R>0$. By \corref{compact}, this implies relative compactness of
	$\set{\smallx_n}{n\in\N}$ in the Gromov-vague topology. Hence the sequence converges Gromov-vaguely.
\end{proof}\sm

Unfortunately, the subspaces $\XHB$ and $\Xc$ are not Polish, and hence it is impossible to find a complete
metric inducing Gromov-vague topology on them. They are, however, Lusin spaces. Recall that a metrizable
topological space is, by definition, a Lusin space if it is the image of a Polish space under a continuous,
bijective map.

\begin{cor}[$\XHB$ and $\Xc$ are Lusin]
	The space\/ $\XHB$ of Heine-Borel locally finite measure spaces, equipped with the Gromov-vague
	topology, is a Lusin space but not Polish. The same is true for the space\/ $\Xc$ of compact metric
	finite measure spaces.
\label{C:Lusin}
\end{cor}\sm

\begin{proof}
	$\XHB$ and $\Xc$ are measurable subsets (\corref{measurable}) of the Polish space $\Xmm$. Hence they are
	Lusin by Theorem 8.2.10 of \cite{Cohn1980}.

	To see that $\XHB$ is not Polish, note that it is a dense subspace of $\Xmm$, and using \lemref{HBchar}
	we see that its complement, $\Xmm\setminus\XHB$, contains a countable intersection of open dense sets, namely
	$G:=\bigcap_{a>0,\,a\in\Q} \bset{\smallx\in\Xmm}{\lmbfct11(\smallx) < a}$.
	Such a subspace cannot be Polish by standard arguments (see also \cite[Remark~4.7]{Loehr}), which we
	recall for the reader's convenience.
	Assume for a contradiction that $\XHB$ is Polish. By the Mazurkiewicz theorem (\cite[Theorem~8.1.4]{Cohn1980}),
	it is a countable intersection of open sets, $\XHB=\bigcap_{n\in\N} U_n$, say. Obviously, the $U_n$ have
	to be dense, because $\XHB$ is. Now $\XHB \cap G$ is also a countable intersection of open dense sets,
	hence it is dense by the Baire category theorem (\cite[Theorem~D.37]{Cohn1980}). This is a contradiction,
	because $G\subseteq \Xmm\setminus \XHB$.

	The same reasoning also applies to $\Xc$, hence $\Xc$ is also not Polish.
\end{proof}\sm

\section{The Gromov-Hausdorff-vague topology}
\label{S:GHvague}
In this section we introduce with the {\em Gromov-Hausdorff-vague topology} a topology which is stronger than the Gromov-vague topology.
The need for such a topology can be motivated by situations as in Example~\ref{ex:conv}, and by the fact that
there are sequences of finite measures $(\mu_n)_{n\in\mathbb{N}}$ on a common compact space $(E,d)$, such that
$\mu_n\Rightarrow\mu$, as $n\to\infty$, but their supports do not converge. The convergence of supports, however, plays
a crucial r\^ole for the convergence of associated random walks to a Brownian motion on the limit space (see
\cite{AthreyaLoehrWinter}).
We define the stronger topology based on isometric embeddings, discuss its connection to the related measured
(Gromov-)Hausdorff topology and to the Gromov-Hausdorff-Prohorov metric known from the literature, state a
stability result, and characterize compact sets. A main result of this section is Polishness of the
Gromov-Hausdorff-weak and -vague topologies (Propositions~\ref{P:GHwPolish} and \ref{Prop:Polish}). Interpreted
in terms of the Gromov-Hausdorff-Prohorov metric as used in \cite{Miermont2009}, this means that the subspace of
metric measure spaces with full support of the measure is Polish although it is not closed (\corref{sGHP}).

We cannot, of course, build such a strong notion of convergence on the notion of sampling alone, and
therefore rather use an isometric embedding approach (compare Proposition~\ref{P:embed}).
Recall that the Hausdorff distance between two closed subsets $A,B$ of a metric space $(E,d)$ of bounded diameter is defined by
\begin{equation}
\label{e:Hausdorff}
  \dH\big(A,B\big)
 :=
   \inf\big\{\eps>0:\,A^\eps\supseteq B,\mbox{ and }B^\eps\supseteq A\big\},
\end{equation}
where once more $A^\eps:=\{x\in A:\,d(x,A)\le\eps\}$ denotes the closed $\eps$-neighborhood of $A$.
Recall that $\Xfin$ denotes the space of metric finite measure spaces.

\begin{definition}[(pointed) Gromov-Hausdorff-weak topology] Let for each\/ $n\in\mathbb{N}\cup\{\infty\}$,
$\tree_n:=(X_n,r_n,\rho_n,\mu_n) \in \Xfin$. We say that\/ $(\tree_n)_{n\in\mathbb{N}}$ converges to $\tree_\infty$ in
\emph{Gromov-Hausdorff weak topology} if and only if there exists a pointed metric space $(E,d_E,\rho_E)$
and, for each $n\in\mathbb{N}\cup\{\infty\}$, an isometry
$\varphi_n\colon \supp(\mu_n)\to E$ with\/
$\varphi_n(\rho_n)=\rho_E$, and such that in addition to
\begin{equation}
\label{e:imageweak}
   (\varphi_n)_\ast\mu_n\Tno(\varphi_\infty)_\ast\mu_\infty,
\end{equation}
also
\begin{equation}
\label{e:Gweak33}
   \dH\Bigl(\vphi_n\(\supp(\mu_n)\),\, \varphi_\infty\(\supp(\mu_\infty)\)\Bigr)\tno 0.
\end{equation}
\label{Def:GHw}
\end{definition}\sm

A very similar topology for compact metric measure spaces was first introduced in \cite{Fukaya1987} under the
name \emph{measured Hausdorff topology} (often referred to as measured Gromov-Hausdorff topology) and further discussed
in \cite{EvansWinter2006,Miermont2009}.
The definition of this topology is exactly the same as that of the Gromov-Hausdorff-weak topology, except that
$\supp(\mu_n)$ is replaced by $X_n$, $n\in\N\cup\{\infty\}$. The consequence is that, when comparing compact metric
measure spaces, the geometric structure outside the support is taken into account, while it is ignored by our
definition. It is important to note that this leads to different equivalence classes, i.e., the measured Hausdorff
topology is \emph{not} defined on $\Xc$, but rather on a space of equivalence classes with respect to the following
equivalence relation.
We say that two metric measure spaces $(X,r,\rho,\mu)$ and $(X',r',\rho',\mu')$ are \define{strongly equivalent}
if and only if there is a surjective isometry $\phi\colon X\to X'$ such that $\phi(\rho)=\rho'$ and
$\phi_\ast(\mu)=\mu'$. Define
\begin{equation}
	\XmHc := \{\,\text{strong  equivalence classes of compact metric measure spaces}\,\}.
\end{equation}
it is well-known that the measured Hausdorff topology is induced by the so-called
\emph{Gromov-Hausdorff-Prohorov metric} defined on $\XmHc$ as follows. For $\smallx=(X,r,\rho,\mu)$,
$\smallx'=(X',r',\rho',\mu') \in \XmHc$,
 \begin{equation}
 \label{e:DGPrcomp}
	\dGHP(\smallx, \smallx') :=
   \inf_{d}\, \dPr[(X\sqcup X',d)](\mu,\mu') +\dH[(X\sqcup X',d)](X,X') +d(\rho,\rho'),
\end{equation}
where the infimum is taken over all metrics $d$ on $X\sqcup X'$ that extend both $r$ and $r'$. Note that $(\XmHc,
\dGHP)$ is a complete, separable metric space (see \cite[Proposition~8]{Miermont2009}).

Now we can easily identify $\Xc$ with the subspace of $\XmHc$ that consists of all (strong equivalence classes of)
compact metric spaces with a measure of full support, i.e.\ with
\begin{equation}
	\XmHcsupp := \bset{(X, r, \rho, \mu) \in \XmHc}{X=\supp(\mu)},
\end{equation}
by choosing representatives with full support from the larger equivalence classes of $\Xc$, i.e.\ via the injective map
\begin{equation}\label{e:iota}
	\iota \colon \Xc \to \XmHcsupp,\qquad (X, r, \rho, \mu) \mapsto \(\supp(\mu), r, \rho, \mu\).
\end{equation}
It is obvious that $\iota$ is a homeomorphism if we equip $\Xc$ with the Gromov-Hausdorff-weak and $\XmHcsupp$ with the
measured Hausdorff topology. Its inverse $\iota^{-1}$ can naturally be extended to all of $\XmHc$, but this
extension looses continuity, as we show in the following remark.

\begin{remark}[support projection]
	Equip $\XmHc$ with the measured Hausdorff topology and $\Xc$ with the Gromov-Hausdorff-weak topology. The \emph{support projection}
\begin{equation}
	\pi^\mathrm{supp} \colon \XmHc \to \XmHcsupp,\qquad (X, r, \rho, \mu) \mapsto \(\supp(\mu), r, \rho, \mu\).
\end{equation}
	is an open map, but neither continuous nor closed. In particular, associating to a strong equivalence
	class of metric measure spaces in $\XmHc$ the corresponding equivalence class in $\Xc$ is not
	a continuous operation, although it induces a homeomorphism from $\XmHcsupp$ onto $\Xc$.
\end{remark}\sm

\begin{remark}[full support assumption]
	The requirement that the measure on a metric space has full support is not unnatural.
	It plays, for instance, a crucial r\^ole for defining Markov processes via Dirichlet forms
	(a particular example is \cite{AthreyaLoehrWinter}), and is even included in the definition of ``Radon measure'' in
	\cite{FukushimaOshimaTakeda1994}.
\label{Rem:004}
\hfill$\qed$
\end{remark}\sm

Note that $\XmHcsupp$ is not closed in $\XmHc$, hence transporting the Gromov-Hausdorff-Prohorov metric $\dGHP$ with
$\iota$ back to $\Xc$ does not yield a complete metric. The following proposition shows, however, that we can find a
different, complete metric for the Gromov-Hausdorff-weak topology on $\Xc$. This also implies that, although
$(\XmHcsupp, \dGHP)$ is not complete, it can still be used as a Polish state-space, because the induced
topological space is Polish. To define the complete metric on $\Xc$, we use the global lower mass function $\glmb$
from \eqref{glmb}, the Gromov-Hausdorff-Prohorov metric $\dGHP$ from \eqref{e:DGPrcomp}, and the homeomorphism
$\iota$ from \eqref{e:iota}. Recall that $\glmb > 0$ on $\Xc$ for every $\delta>0$ by \lemref{cchar}.

\begin{definition}
	For\/ $\smallx, \smallx' \in \Xc$, let
	\begin{equation}
		\dsGHP(\smallx, \smallx') := \dGHP\(\iota(\smallx), \iota(\smallx')\) + \integral01{1\land
			\left|\frac1{\glmb(\smallx)} - \frac1{\glmb(\smallx')}\right|}\delta.
	\end{equation}
	We call\/ $\dsGHP$ the \define{support Gromov-Hausdorff-Prohorov metric}.
\end{definition}\sm

\begin{proposition}[$(\Xc,\dsGHP)$ is a complete metric space]
	The metric $\dsGHP$ induces the Gromov-Hausdorff-weak topology on\/ $\Xc$.
	Furthermore, $(\Xc, \dsGHP)$ is a complete, separable metric space.
\label{P:GHwPolish}
\end{proposition}\sm

\begin{proof}
	Let $(\smallx_n)_{n\in\mathbb{N}}$ and $\smallx$ be in $\Xc$. Then $\glmb(\smallx) >0$,
for all $\delta>0$. Thus by definition, $\dsGHP(\smallx_n, \smallx) \tno 0$ if and only if
	\begin{equation}
		\dGHP\(\iota(\smallx_n),\iota(\smallx)\) \tno 0,
\label{eq:iotaconv}
\end{equation}
and for almost all $\delta>0$,
\begin{equation}
		 \glmb(\smallx_n) \tno \glmb(\smallx).
\label{eq:lmbconv}
\end{equation}
	Because $\iota$ is a homeomorphism, \eqref{eq:iotaconv} is equivalent to the Gromov-Hausdorff-weak
	convergence $\smallx_n \tno \smallx$. We have to show that this already implies \eqref{eq:lmbconv}, i.e.\
	that $\smallx$ is continuity point of $\glmb$ w.r.t.\ Gromov-Hausdorff-weak topology for almost all $\delta>0$.
	To see this, recall that $\glmb$ is upper semi-continuous w.r.t.\ Gromov-vague topology
	(\lemref{lmbusc}), and a fortiori also w.r.t.\ Gromov-Hausdorff-weak topology. Assume that all
	$\smallx_n=(X_n, r_n, \rho_n, \mu_n)$ and $\smallx=(X,r,\rho, \mu)$ are embedded in some common space
	$(E, d, \rho)$ such that $\mu_n$ converges weakly to $\mu$ and $\supp(\mu_n)$ in Hausdorff metric to
	$\supp(\mu)$. Then, for every $\deltah<\delta$ and $n$ sufficiently large, every $\delta$-ball around
	some $y\in \supp(\mu_n)$ contains a $\deltah$-ball around some $x\in \supp(\mu)$. Therefore, $\nliminf
	\glmb(\smallx_n) \ge \glmb[\deltah](\smallx)$. This means that $\glmb$ is Gromov-Hausdorff-weakly lower
	semi-continuous in $\smallx$ for every $\delta>0$ with $\glmb(\smallx) = \sup_{\deltah < \delta}
	\glmb[\deltah](\smallx)$. Because $\delta \mapsto \glmb(\smallx)$ is an increasing function, this is the
	case for almost all $\delta>0$. This means that \eqref{eq:lmbconv} is implied by Gromov-Hausdorff-weak
	convergence, and hence $\dsGHP$ induces Gromov-Hausdorff-weak topology as claimed.\sm

	That $(\Xc, \dsGHP)$ is a {\em separable} metric space is obvious, and it remains to show its {\em completeness}.
Consider a $\dsGHP$-Cauchy sequence
	$(\smallx_n)_{n\in\N}$ in $\Xc$. Then, by completeness of $\dGHP$ on $\XmHc$,
	the sequence $\(\iota(\smallx_n)\)_{n\in\N}$ converges in measured Hausdorff topology to some
	$\smally=(X, r, \rho, \mu)\in \XmHc$. We have to show $\smally \in \iota(\Xc)=\XmHcsupp$.
	Assume for a contradiction that this is not the case, i.e.\ there exists $x\in X \setminus \supp(\mu)$.
	Then there is a $\delta>0$ with $B(x, 2\delta)\cap \supp(\mu) = \emptyset$. By the measured Hausdorff
	convergence and the fact that $\iota(\smallx_n) \in \XmHcsupp$ for all $n$, this clearly implies
	$\glmb(\smallx_n) \tno 0$. This, however, cannot be the case because $(\smallx_n)_{n\in\mathbb{N}}$ is a Cauchy sequence
	w.r.t.\ $\dsGHP$.
\end{proof}\sm

\begin{cor}\label{C:sGHP}
	The set $\XmHcsupp$ of (strong equivalence classes of) compact metric full-support measure spaces with
	the topology induced by the Gromov-Hausdorff-Prohorov metric $\dGHP$ is a Polish space (although $\dGHP$
	restricted to $\XmHcsupp$ is not complete).
\end{cor}\sm

\begin{cor}[Gromov-Hausdorff-weak compactness]\label{C:GHWcompact}
	A set\/ $\K\subseteq \Xc$ is relatively compact in the Gromov-Hausdorff-weak topology if and only if the
	following hold
	\begin{enumerate}
		\item\label{i:mb} The set of the total masses is uniformly bounded, i.e.,
\begin{equation}
\label{e:totalunif}\sup_{(X,r,\rho, \mu)\in \Xc} \mu(X)< \infty.
\end{equation}
		\item\label{i:covnum} For all\/ $\eps>0$ there exists an $N_{R,\eps}\in\N$ such that for
			all\/ $(X,r,\rho,\mu)\in\K$, $\supp(\mu)$ can be covered by $N_{R,\eps}$ many balls of
			radius $\eps$.
		\item\label{i:mass} $\K$ satisfies the global lower mass-bound property.
	\end{enumerate}
\end{cor}\sm

\begin{proof}
	From \propref{GHwPolish}, the definition of $\dsGHP$, and the fact that $\dGHP$ induces the measured
	Hausdorff topology, we see that $\K$ is Gromov-Hausdorff-weakly relatively compact if and only if
	$\iota(\K)$ is relatively compact in measured Hausdorff topology, and $1/\glmb$ is bounded on $\K$. The
	latter is obviously equivalent to the global lower mass-bound \ref{i:mass} If $\K\subseteq \Xprob$, the measured
	Hausdorff relative compactness of $\iota(\K)$ is equivalent to \ref{i:covnum} by
	\cite[Proposition~2.4]{EvansWinter2006} together with \cite[Theorem~7.4.15]{BurBurIva01}.
	It is therefore easy to see that it is in general equivalent to \ref{i:covnum} together with \ref{i:mb} (compare
	\cite[Remark~7.2(ii)]{GrevenPfaffelhuberWinter2009}).
\end{proof}

In the same way as we used the Gromov-weak topology to define the Gromov-vague topology, we also define the
Gromov-Hausdorff-vague topology on $\Xmm$ based on the Gromov-Hausdorff-weak topology on $\Xfin$.
\begin{definition}[(pointed) Gromov-Hausdorff-vague topology]
	Let for each\/ $n\in\mathbb{N}\cup\{\infty\}$, $\tree_n:=(X_n,r_n,\rho_n,\mu_n)$ be in $\mathbb{X}$.
	We say that\/ $(\tree_n)_{n\in\mathbb{N}}$ converges to $\tree_\infty$ in \emph{Gromov-Hausdorff
	vague topology} if and only if\/ $\treenR\ton(\smallx_\infty)\restricted{R}$ Gromov-Hausdorff-weakly for all but
	countably many $R>0$.
\label{Def:002}
\end{definition}\sm

The following embedding result and its corollary about Gromov-Hausdorff-vaguely compact sets are proved in the same way
as \propref{embed} and \corref{compact}.
\begin{proposition}[isometric embeddings; Gromov-Hausdorff-Prohorov metric] Let for each $n\in\mathbb{N}\cup\{\infty\}$,
$\tree_n:=(X_n,r_n,\rho_n,\mu_n)$ be
in $\XHB$. The following are equivalent:
\begin{enumerate}
\item
$\tree_n \tno\tree_\infty$,
Gromov-Hausdorff vaguely.
\item There exists a rooted Heine-Borel
space $(E,d_E,\rho_E)$ and for each $n\in\mathbb{N}\cup\{\infty\}$ isometries $\varphi_n\colon \supp(\mu_n)\to E$ with
$\varphi_n(\rho_n)=\rho_E$, and such that in addition to \eqref{e:embed}, also
\begin{equation}
\label{e:Gweak3}
   d_{\mathrm{H}}\Big(\varphi_n(\supp\mu_n)\cap \Bcl_{d_E}(\rho_E, R),\,\vphi_\infty(\supp\mu_\infty)\cap \Bcl_{d_E}(\rho_E, R)\Big)\tno 0,
\end{equation}
for all but countably many $R>0$.
\item $\dasGHP(\smallx_n, \smallx_\infty)\tno 0$, where for $\smallx,\smallx'\in\XHB$,
\begin{equation}
\label{e:002}
	\dasGHP\big(\smallx, \smallx'\big) :=
		\plainint{e^{-R}\(1\land \dsGHP(\smallx\restricted{R}, \smallx'\restricted{R})\)}R.
\end{equation}
\end{enumerate}
\label{P:embedGHv}
\end{proposition}\sm

\begin{cor}[Gromov-Hausdorff-vague compactness]
For a set $\K\subseteq \Xmm$ the following are equivalent:
	\begin{enumerate}
	\item $\K$ is relatively compact in $\Xmm$ equipped with the Gromov-Hausdorff-vague topology.
	\item For all\/ $R>0$, the set of restrictions $\K\restricted{R} :=
		\bset{\treeR}{\tree\in \K}$ is relatively compact in the Gromov-Hausdorff-weak topology.
	\item $\K\restricted{R_k}$ is relatively compact in the Gromov-Hausdorff-weak topology for a
		sequence $R_k\to\infty$.
	\end{enumerate}
\label{C:005}
\end{cor}\sm

\begin{remark}[Gromov-Hausdorff-Prohorov and length spaces]
Under the name \emph{Gromov-Hausdorff-Prohorov topology}, the measured Hausdorff topology was recently extended in
\cite{AbrahamDelmasHoscheit2013} to the space of complete, locally compact length spaces equipped with locally finite
measures. The extension was done with the same localization procedure that we use. Note the following:
\begin{enumerate}
\item Complete locally compact length spaces are Heine-Borel spaces and well suited for applications concerning $\R$-trees.
	The assumption of being a length space and thereby path-connected, however, is too restrictive in general.
	For example, in Theorem~1 of \cite{AthreyaLoehrWinter} we establish convergence in path space of continuous time
	random walks on discrete trees to time-changed Brownian motion on $\mathbb{R}$-trees (appearing as the
	Gromov-Hausdorff-vague limit of the discrete trees), where the underlying trees are encoded as metric spaces and
	jump rates and/or time-changes are encoded by the so-called speed measure. Since we need the speed measure to
	have full support, the situation is incompatible with a connectedness requirement.
\item In a general setting, the name Gromov-Hausdorff-Prohorov topology might be a bit misleading, as ``Prohorov''
	suggests weak convergence, while the localized convergence is vague in the sense that mass can get lost. Also
	note that, if we drop the assumption of being length spaces, the localized convergence is not really an
	extension of measured Hausdorff convergence any more (compare Remark~\ref{Rem:002}).
\label{Rem:013}
\hfill$\qed$
\end{enumerate}
\end{remark}\sm

\comment{
For characterizing compact sets, one can rely on a description of compact sets in the Gromov-Hausdorff topology and then insert an extra condition on the measure.
\begin{remark}[Gromov-Hausdorff-vague compactness]
	Compactness criteria for the Gromov-weak topology on metric probability measure spaces are given in Proposition~2.4
	in \cite{EvansWinter2006}. Combining this with Theorem~7.4.15 in \cite{BurBurIva01} and Remark~7.2(ii)
	in \cite{GrevenPfaffelhuberWinter2009} yields that\/ $\K$ is relatively compact in $\Xmm$ equipped with
	the Gromov-Hausdorff-vague topology if and only if for all\/ $R>0$ the following hold:
\begin{enumerate}
\item The total mass of $B(\rho,R)$ is uniformly bounded, i.e., (\ref{e:umb}) holds.
\item For all\/ $\varepsilon>0$ there exists a $N_{R,\varepsilon}\in\mathbb{N}$ such that for all\/ $(X,r,\rho,\mu)\in\K$,
	$\supp(\mu) \cap B(\rho,R)$ can be covered by $N_{R,\varepsilon}$ many balls of radius $\varepsilon$.
\end{enumerate}

In particular, if $\bar{\K}\subseteq \XHB$, then $\K$ is relatively compact in the  Gromov-Hausdorff-vague topology
if and only if  $\K$ is relatively compact in the  Gromov-vague topology.
\label{Rem:005}
\hfill$\qed$
\end{remark}\sm
}

\begin{proposition}[$\XHB$ with Gromov-Hausdorff-vague topology is Polish]
	The space $\XHB$ of Heine-Borel boundedly finite measure spaces equipped with the Gromov-Hausdorff-vague topology is a Polish space.
\label{Prop:Polish}
\end{proposition}\sm

\begin{proof}
	We follow the proof of \propref{Polish} and define
	\begin{equation}
		\dasGHP\big(\smallx,\smally\big)
 :=
		\inta{\R_+}{e^{-R}\(1\land \dsGHP(\smallx\restricted{R}, \smally\restricted{R})\)}R.
	\end{equation}
	We know from \propref{embedGHv} that $\dasGHP$ induces the Gromov-Hausdorff-vague topology.
	Separability and completeness follow from the corresponding properties of $\dsGHP$ (\propref{GHwPolish}) and the
	compactness criterion given in \corref{005}, in the same way as in the proof of \propref{Polish}.
\end{proof}\sm

Even though the Gromov-Hausdorff-vague topology is nice (i.e.\ Polish) on $\XHB$ and defined on all of $\Xmm$, it
appears to be too strong to be useful on the larger space.
\begin{remark}[Gromov-Hausdorff-vague topology is non-separable on\/ $\Xmm$]
The spaces	$\Xmm$ and\/ $\Xprob$, equipped with the Gromov-Hausdorff-vague topology, are not separable. In particular they are
	not Lusin spaces.
	Indeed, we can topologically embed the non-separable space $l^\infty_+$ into $\Xprob$ as follows:
	for $n\in\N$ and $a\in \R_+$, let\/ $A_a^n:=\{n\}\times [0,a]^n$, and $\mu_a^n$ some measure on $A_a^n$ with
	full support and total mass $2^{-n}$.
	Define $\psi\colon l^\infty_+ \to \Xprob$ by
	$\psi(a) := \(\bigcup_{n\in\N} A_{a_n}^n, r, \rho, \sum_{n\in\N}\mu_{a_n}^n\)$,
	where $\rho=(1,0)$, and $r$ is the supremum of the discrete metric on the first component and the Euclidean
	metric on the second component. It is straightforward to check that $\psi$ is a homeomorphism onto its
	image.
\label{Rem:nonsep}
\hfill$\qed$
\end{remark}\sm

We know from Propositions~\ref{P:Polish}, \ref{P:GHwPolish}, and \ref{Prop:Polish} that $\Xmm$ with Gromov-vague topology,
$\Xc$ with Gromov-Hausdorff-weak topology, and $\XHB$ with Gromov-Hausdorff-vague topology are Polish spaces.
Furthermore, it is known from \cite[Theorem~1]{GrevenPfaffelhuberWinter2009} that $\Xprob$ and $\Xfin$ with Gromov-weak topology are Polish. In the case of $\Xprob$, this is also true for the Gromov-vague topology,
although $\Xprob$ is not Gromov-vaguely closed in $\Xmm$ (see Remark~\ref{Rem:002}).
On the other hand, \corref{Lusin} proves that $\Xc$ and $\XHB$ with Gromov-vague topology are Lusin but
\emph{not} Polish. Similar arguments also show that $\Xfin$ with Gromov-vague topology and
$\Xc$ with Gromov-weak as well as with Gromov-Hausdorff-vague topology are Lusin and not Polish.
Gromov-Hausdorff-vague topology is not even separable on $\Xmm$ and $\Xprob$ by Remark~\ref{Rem:nonsep}. We
summarize the situation in Figure~\ref{fig:Poltable}.

\begin{figure}
\begin{tabular}{c|cccccc}
	    & $\Xmm$    & $\XHB$ & $\Xc$ & $\Xfin$    & $\Xprob$  &$\Xc\cap \Xprob$\\
	\hline
	Gv  & Polish    & Lusin  & Lusin  & Lusin     & Polish    & Lusin \\
	Gw  & --        & --     & Lusin  & Polish    & Polish    & Lusin \\
	GHv & non-sep.\ & Polish & Lusin  & non-sep.\ & non-sep.\ & Lusin \\
	GHw & --        &  --    & Polish & non-sep.\ & non-sep.\ & Polish
\end{tabular}
\caption{The table shows topological properties of different spaces of metric measure spaces in different
topologies. Entries: ``--'' means not defined; ``non-sep.'' means non-separable; ``Lusin'' means Lusin but \emph{not} Polish.
Spaces: the spaces are defined in \defref{001} and Remark~\ref{Rem:002}.
Topologies: Gv=Gromov-vague, Gw=Gromov-weak, GHv=Gromov-Hausdorff-vague, GHw=Gromov-Hausdorff-weak.}
\label{fig:Poltable}
\end{figure}

We conclude this section with the following stability property, which is the analogue of Lemma~\ref{Lem:001}. It is an
immediate consequence of the definition of the Gromov-Hausdorff-vague topology by means of isometric embeddings. The
proof follows the same lines as the proof of Lemma~\ref{Lem:001} and is therefore omitted.
\begin{lemma}[perturbation of measures]
	Consider\/ $ \smallx = (X, r, \rho, \mu),\,\smallx_n = (X_n, r_n, \rho_n, \mu_n) \in \Xmm$, and
	another boundedly finite measure $\mu_n'$ on $X_n$, $n\in\N$. Assume that $\smallx_n \ton \smallx$
	Gromov-Hausdorff-vaguely, and that there exists a sequence $R_k\to \infty$ such that for all\/ $k\in\mathbb{N}$,
	\begin{equation}\label{Prmerg2}
		d_{\mathrm{Pr}}^{(X_n,r_n)}\bigl(\mu_n\restricted{R_k}, \mu_n'\restricted{R_k}\bigr)\tno 0,
			\>\mbox{ and\/ }\>
		d_{\mathrm{H}}^{(X_n,r_n)}\bigl(\supp(\mu_n\restricted{R_k}),\supp(\mu_n^\prime\restricted{R_k})\bigr)\tno 0.
	\end{equation}
	Then $(X_n, r_n, \rho_n, \mu_n')$ converges Gromov-Hausdorff-vaguely to $\smallx$.
\label{Lem:pertGHvag}
\end{lemma}\sm

\begin{example}[Normalized length measure versus degree measure]
	Consider a graph theoretic tree\/ $T'$ which is locally finite, i.e.\  $\deg(v) < \infty$ for all\/ $v\in
	T'$, where $\deg$ is the degree of a node. Equip $T'$ with the graph distance $r'$, i.e.\ the length of
	the shortest path, and fix a root\/ $\rho'\in  T'$. Recall the notion of\/ $\R$-tree from
	Example~\ref{Exp:001}. It is well known that\/ $(T', r')$ can be embedded isometrically
	into a complete, locally compact\/ \Rtree\ $(T, r)$ in an essentially unique way.
	Denote the image of\/ $\rho'$ by $\rho$ and the image of\/ $T'$ by $\nod(T)$.
	On $T'$, we consider two natural measures. The \emph{node measure} $\munodT[T']$, which is just the counting
	measure on the nodes (except the root), and the \emph{degree measure} $\mudegT[T']$, which is proportional to
	the degree of the node.  The push-forwards on $T$ are given by
	\begin{equation}
		\munodT := \sum_{x\in \nod(T)\setminus\{\rho\}} \delta_x
			\und
		\mudegT := \tfrac12 \sum_{x\in \nod(T)} \deg(x) \cdot \delta_x.
	\end{equation}
	Note that $(T', r', \rho', \mudegT[T']) \cong (T, r, \rho, \mudegT)$, and similarly for the node measure.
	On $T$, there is also a third natural measure, namely the \emph{length measure} $\lambda=\lambdaT$,
	which is the 1-dimensional Hausdorff measure on $T\setminus \lf(T)$, where\/ $\lf(T)=\set{x\in
	T}{T\setminus\{x\} \text{ is connected}}$ is the set of leaves of\/ $T$. Note that
	$\lambdaT(T)=\munodT(T)=\mudegT(T)$.

	Now consider a sequence $(T'_n)_{n\in\N}$ of locally finite, graph theoretic trees, and the
	rooted \Rtree s $(T_n, r_n, \rho_n)$ constructed as above. We assume that there are two sequences
	$(\alpha_n)_{n\in\N}$ and\/ $(\beta_n)_{n\in\N}$ of positive numbers, both of which converge to $0$,
	such that
	\begin{equation}\label{eq:lengthconv}
		\(T_n, \alpha_n r_n, \rho_n, \beta_n \lambdaTn\) \tno \smallx,
	\end{equation}
	Gromov-Hausdorff-vaguely for some $\smallx=(T,r,\rho,\mu)\in\XHB$, which is necessarily an \Rtree.
	Such a convergence can often be deduced via convergence of excursions, see \propref{glue} and
	Example~\ref{Exp:003} below. We claim that in this case, the length measure can be replaced by the degree
	measure or the node measure, i.e.\ that \eqref{eq:lengthconv} implies the Gromov-Hausdorff-vague convergences
	\begin{equation}\label{eq:replace}
		\(T_n, \alpha_n r_n, \rho_n, \beta_n \mudegTn\) \tno \smallx
			\und
		\(T_n, \alpha_n r_n, \rho_n, \beta_n \munodTn\) \tno \smallx.
	\end{equation}
	
	Indeed, we have\/ $\supp(\munodTn) = \supp(\mudegTn) = \nod(T_n)$, $\supp(\lambdaTn)=T_n$ and, for every $R>0$,
	\begin{equation}\label{eq:H}
		\dH\(\nod(T_n)\cap B(\rho_n, R),\, B(\rho_n, R)\)
		\le \alpha_n \tno 0.
	\end{equation}
	For the Prohorov distance, assume first that the diameter of\/ $T_n$ is smaller than $R$. Then
	\begin{equation}
		\dPrTn(\mudegTn, \lambdaTn) \le \half\alpha_n \und
		\dPrTn(\munodTn, \lambdaTn) \le \alpha_n.
	\end{equation}
	In the general case, we have to take boundary effects into account. Using the annulus
	$S^\eps(\rho_n,R) := \Bcl(\rho_n, R+\half\eps)\setminus B(\rho_n, R-\half\eps)$, we obtain
	\begin{equation}\label{eq:Pr}
		\dPrTn(\mudegTn\restricted R,\, \lambdaTn\restricted R) \,\le\, \half\alpha_n \lor
			\beta_n\cdot\lambdaTn\(S^{\alpha_n}(\rho_n, R)\),
	\end{equation}
	and a similar estimate for\/ $\munodTn$ instead of\/ $\mudegTn$.
	Using\/ \eqref{eq:lengthconv} we see that\/ $\beta_n\lambdaTn\(S^{\alpha_n}(\rho_n, R)\)$ tends to zero for all\/
	$R$ with\/ $\mu\(S(\rho,R)\)=0$. Therefore the claimed Gromov-Hausdorff-vague convergences\/ \eqref{eq:replace}
	follow from \eqref{eq:Pr}, \eqref{eq:H} and Lemma~\ref{Lem:pertGHvag}.
\label{Exp:002}
\hfill$\qed$
\end{example}\sm

\section{Closing the gap}
\label{S:relation}
In this section we prove the main criterion for convergence in Gromov-Hausdorff-vague topology.
We shall use notation used in (\ref{d:002a}), (\ref{d:002}) and the definitions of the lower mass functions $\lmb$ and $\glmb$ from \eqref{lmb} and \eqref{glmb}, respectively.

In order for a sequence $\(\smallx_n\)_{n \in \N}:=\(X_n,r_n,\rho_n,\mu_n\)_{n\in\N}$ of compact metric finite measure spaces to
converge in Gromov-Hausdorff-weak topology to a space $\smallx=(X,r,\rho,\mu) \in \Xc$, it certainly has to converge in
the weaker Gromov-weak topology. This is a kind of ``finite-dimensional convergence'', which is expressible in terms of
sampling finite sub-spaces:
\begin{enumerate}
	\item\label{it:fd}  For all $k\in\mathbb{N}$, and $\vphi\in\bar{{\mathcal C}}(\R_+^{k+1\choose 2})$
	\begin{equation}\label{e:003}
	\begin{aligned}
	   &\int\mu_n^{\otimes k}(\mathrm{d}(x^n_1,\ldots,x^n_k))\,\vphi\big((r_n(x^n_i,x^n_j))_{0\le i<j\le k}\big)
	   \\
	   &\tno
	   \int\mu^{\otimes k}(\mathrm{d}(x_1,\ldots,x_k))\,\vphi\big((r(x_i,x_j))_{0\le i<j\le k}\big)
	\end{aligned}
	\end{equation}
	where we put $x^n_0:=\rho$ and $x_0:=\rho$.
\end{enumerate}
We show in \thmref{mGH} below that, given \ref{it:fd}, Gromov-Hausdorff-weak convergence follows from a simple
``tightness condition'', which is given in terms of the lower mass function:
\begin{enumerate}\setcounter{enumi}{1}
	\item\label{it:tight} For all $\delta>0$, $\liminf_{n\to\infty}\glmb(\smallx_n)>0$.
\end{enumerate}
Note that for checking \ref{it:fd} and \ref{it:tight}, we do not have to find any embedding into a common metric space.
We actually show that \ref{it:fd} and \ref{it:tight} together are even equivalent to Gromov-Hausdorff-weak convergence,
and this characterization even holds if the $\smallx_n$ are not compact (but $\smallx$ is).

\begin{theorem}[Gromov-weak versus Gromov-Hausdorff-weak convergence]
	Let\/ $\smallx = (X, r, \rho, \mu)$ and\/ $\smallx_n = (X_n,r_n,\rho_n,\mu_n)$, $n\in\N$, be
	metric finite measure spaces. Then the following are equivalent.
	\begin{enumerate}
	\item\label{i:Gw} $(\smallx_n)_{n\in\N}$ converges in Gromov-weak topology to\/ $\smallx$, and for all $\delta>0$,
		\begin{equation} \label{e:udensity}
			\nliminf\glmb(\smallx_n)>0.
		\end{equation}
	\item\label{i:mGH} $\smallx$ is compact, and\/ $(\smallx_n)_{n\in\N}$ converges in
		Gromov-Hausdorff-weak topology to\/ $\smallx$.
	\end{enumerate}
	If\/ $\smallx_n$ is compact for all\/ $n\in\N$, the following is also equivalent:
	\begin{enumerate}\setcounter{enumi}{2}
	\item\label{i:Gwlm} $(\smallx_n)_{n\in\N}$ converges in Gromov-weak topology to\/ $\smallx$, and\/
		$\set{\smallx_n}{n\in\N}$ satisfies the global lower mass-bound property
		(Definition~\ref{Deef:001}).
	\end{enumerate}
\label{T:mGH}
\end{theorem}\sm

\begin{proof} \emph{``\ref{i:mGH}$\Rightarrow$\ref{i:Gw}'' }
	Assume $(\supp(\mu_n),r_n,\rho_n,\mu_n)\tno(\supp(\mu),r,\rho,\mu)$ Gromov-Hausdorff-weakly.
	W.l.o.g.\ we may assume that $X=\supp(\mu)$, $X_n=\supp(\mu_n)$, and that $X_n$ and $X$ are embedded into a
	complete, separable metric space $(E, d)$ such that
	\begin{equation} \label{e:005}	
	   d^{(E,d)}_{\mathrm{Pr}}\big(\mu_n,\mu\big)\tno 0, \quad\mbox{and}\quad d^{(E,d)}_{\mathrm{H}}\big(X_n,X\big) \tno 0.
	\end{equation}
	Furthermore let $(X,r)$ be compact. We need to show (\ref{e:udensity}).
	Assume to the contrary  that there exists $\delta>0$ and
	$x_n\in X_n$ such that $\nliminf \mu_n\(B(x_n, 2\delta)\)=0$.
	Due to (\ref{e:005}) we can find $y_n \in X$
	with $d(x_n, y_n)\tno 0$. Moreover by (\ref{e:005}),
	\begin{equation} \label{e:004}
		\nliminf\mu\big(B(y_n, \delta)\big) \le \nliminf\mu_n\big(B(x_n, 2\delta)\big) = 0.
	\end{equation}
	As $X$ is compact, we may assume w.l.o.g.\ that $y_n$ converges
	to some $y\in X$.
	Then $\mu\(B(y, \delta)\) \le \liminf\mu\(B(y_n, \delta)\) = 0$, which contradicts $X=\supp(\mu)$.\sm

\emph{``\ref{i:Gw}$\Rightarrow$\ref{i:mGH}'' } 	Assume that $\smallx_n\tno\smallx$ Gromov-weakly, and
	w.l.o.g.\ that $X=\supp(\mu)$, $X_n=\supp(\mu_n)$, and that $X_n$ and $X$ are embedded into a
	complete, separable metric space $(E, d)$ such that
	$d^{(E,d)}_{\mathrm{Pr}}\big(\mu_n,\mu\big)\tno 0$,
	and that (\ref{e:udensity}) holds. Then, for all $\eps>0$, we can find $n_0=n_0(\eps)\in\mathbb{N}$ such
	that for all $n\ge n_0$,
	\begin{equation}
	\label{e:006}
	   d^{(E,d)}_{\mathrm{Pr}}\big(\mu_n,\mu\big)<\eps \land \inf_{y\in X_n}\mu_n\big(\Bcl(y,\eps)\big)
	   	\land \inf_{x\in X}\mu\(\Bcl(x,\eps)\),
	\end{equation}
	where we used that
	\begin{equation}\label{e:limitlmb}
		\inf_{x\in X}\mu\(\Bcl(x, \eps)\) \ge \nliminf \inf_{y\in X_n} \mu_n\big(\Bcl(y, \tfrac12\eps)\big)>0.
	\end{equation}
	Then, for all $y\in X_n$, $\Bcl(y,\eps)\cap X^\eps\ne\emptyset$, and thus also $\Bcl(y, 2\eps)\cap X \ne
	\emptyset$. Similarly, $\Bcl(x, 2\eps)\cap X_n \ne \emptyset$ for all $x\in X$, and hence
	$d^{(E,d)}_{\mathrm{H}}(X_n,X)\le 2\eps$.

	Compactness of $\smallx$ follows directly from \eqref{e:limitlmb} and \lemref{cchar}.\sm

\emph{``\ref{i:Gwlm}$\Leftrightarrow$\ref{i:Gw}'' }
	Obviously, the global lower mass-bound is equivalent to \eqref{e:udensity} together with
	$\glmb(\smallx_n) > 0$ for all $n\in\N$ and $\delta>0$. The last condition is satisfied for compact
	spaces by \lemref{cchar}.
\end{proof}\sm

The following corollaries are now obvious.

\begin{cor}[Gromov-vague versus Gromov-Hausdorff-vague convergence]
	Let\/ $\smallx = (X, r, \rho, \mu)$ and\/ $\smallx_n = (X_n,r_n,\rho_n,\mu_n)$, $n\in\N$, be
	metric boundedly finite measure spaces. Then the following are equivalent.
	\begin{enumerate}
	\item $(\smallx_n)_{n\in\N}$ converges in Gromov-vague topology to\/ $\smallx$, and for all $\delta>0$ and $R>0$,
		\begin{equation}
			\nliminf\lmb(\smallx_n)>0.
		\end{equation}
	\item $\supp(\mu)$ is Heine-Borel, and\/ $\(\supp(\mu_n),r_n,\mu_n\)_{n\in\N}$ converges in
		Gromov-Hausdorff-vague topology to\/ $\(\supp(\mu),r,\mu\)$.
	\end{enumerate}
	If\/ $\smallx_n$ is Heine-Borel for all $n\in\N$, the following is also equivalent:
	\begin{enumerate}\setcounter{enumi}{2}
	\item $(\smallx_n)_{n\in\N}$ converges in Gromov-vague topology to\/ $\smallx$, and\/
		$\set{\smallx_n}{n\in\N}$ satisfies the local lower mass-bound property (\ref{e:lmb}).
	\end{enumerate}
\label{C:mGH}
\end{cor}\sm

\begin{cor}[Polish subspaces]
	Let\/ $\K\subseteq \XHB$ be a space of Heine-Borel locally finite measure spaces satisfying the local
	lower mass-bound property (\ref{e:lmb}).
	Then its closure\/ $\closure{\K}$ in\/ $\Xmm$ (w.r.t.\ the Gromov-vague topology) is a Polish subspace of\/
	$\XHB$. Furthermore, the Gromov-vague topology and the Gromov-Hausdorff-vague topology coincide on\/
	$\closure{\K}$.
\end{cor}\sm

\begin{cor}[topologies agree up to exceptional sets] Let $\smallx=(X, r, \rho, \mu)$ and $\smallx_n = (X_n,r_n,\rho_n,\mu_n)$, $n\in\N$, be in $\mathbb{X}$.
	Then the following are equivalent:
	\begin{enumerate}
	\item $\smallx_n \tno \smallx$, Gromov-vaguely.
	\item For each $n\in\mathbb{N}$ there is $A_n\subseteq X_n$ such that $\mu_n(A_n \cap B_{r_n}(\rho_n, R))\tno 0$ for all $R>0$, and
\begin{equation}
\label{e:except}		
   \big(X_n\setminus A_n,r_n,\rho_n,\mu_n\restricted{X_n\setminus A_n}\big)\tno\smallx,
\end{equation}
		Gromov-Hausdorff-vaguely.
	\item For each $n\in\mathbb{N}$ there is $A_n\subseteq X_n$ such that $\mu_n(A_n \cap B_{r_n}(\rho_n, R))\tno 0$ for all $R>0$, and (\ref{e:except}) holds
		Gromov-vaguely.
\end{enumerate}
\label{C:except}
\end{cor}\sm


\section{Application to trees coded by excursions}
\label{S:excursions}
In this section we consider encodings of trees by means of excursions.
To be in a position to consider locally compact rather than just compact $\R$-trees we
consider possibly transient excursions,  and conclude from uniform convergence on compacta of a sequence of excursions that the corresponding rooted
boundedly finite $\R$-trees converge Gromov-Hausdorff-vaguely (Proposition~\ref{P:glue}). As an example we present a representation of the scaling limit of a size-biased
Galton-Watson tree (Example~\ref{Exp:003}).

\begin{definition}[(transient) excursions]
A continuous function $e\colon \R_+\to\R_+$ is called \emph{(continuous) excursion} if\/ $e(0)=0$ and $e$ is not
identically $0$. We refer to $\zeta(e):=\sup\bset{s>0}{e(s)>0}$ as the \emph{excursion length}, and to $I_e :=[0,
\zeta(e))$ as the \emph{excursion interval}.
If the excursion length is finite, we call the excursion \emph{compactly supported}. If\/
$\lim_{s\to\infty}e(s)=\infty$, the excursion is called \emph{transient}.
\label{Lem:002}
\end{definition}\sm

Let
\begin{equation}
\label{e:exp}
   \E:=\bigl\{e\colon \R_+\to\R_+ \bigm| e\mbox{ is a continuous excursion}\bigr\}.
\end{equation}
Given $e\in\E$, we define the pseudo-metric $r'_e$ by letting for all $0\le s\le t<\zeta(e)$,
\begin{equation}
   r_e'(s,t):=e(s)+e(t)-2\inf\nolimits_{u\in [s,t]}e(u).
\end{equation}
We write $s\sim_e t$  if $ r_e'(s,t)=0$. Obviously, $s\sim_e t$ is an equivalence relation.

\begin{definition}[glue map]
The \emph{glue map} $g\colon \E\to \Xmm$ sends an excursion to the complete, separable, rooted measure $\R$-tree
\begin{equation}
\label{e:glue}
   g(e) := (T_e,r_e,\rho_e,\mu_e),
\end{equation}
where $T_e := I_e/{\sim_e}$, and $r_e$, $\mu_e$, $\rho_e$ are the push-forwards of $r_e'$, the Lebesgue measure
$\lambda_{I_e}$, and $0$, respectively, under the canonical projection $\pi_e\colon I_e\to T_e$.
\label{Def:glue}
\end{definition}\sm

\begin{lemma}[excursions and associated $\R$-trees] Let $e\in{\mathcal E}$.
\begin{enumerate}
\item If\/ $e$ is compactly supported, then $g(e)$ is a pointed compact finite measure $\R$-tree, i.e.\ $g(e)\in \Xc$.
\item If\/ $e$ is transient, then $g(e)$ is a pointed Heine-Borel boundedly finite measure $\R$-tree, i.e.\ $g(e) \in \XHB$.
\item If\/ $e$ is neither compactly supported nor transient, then $g(e)\not\in\XHB$.
\end{enumerate}
\label{Lem:003}
\end{lemma}\sm

\begin{enproof} \item Follows from Lemma~3.1 in \cite{EvansWinter2006}. \sm

\item Assume that $e$ is transient. Then for all $R>0$, $\xi_e(R):=\sup\{s\ge 0:\,e(s)<R\}<\infty$, and
	$A_R:=\bset{s\in[0,\infty)}{e(s)\le R}$ is a closed subset of $[0,\xi_e(R)]$, and hence compact.
	Note that continuity of $e$ implies continuity of the projection $\pi_e$.
	Therefore, $\Bcl(\rho,R)=\pi_e(A_R)$ is also compact.
	Moreover, $\mu_e\big(\Bcl(\rho,R)\big)\le\xi(R)<\infty$.  As any closed and bounded subset of
	$(T_e,r_e)$ is a closed subset of a closed ball $\Bcl(\rho,R)$ for some $R>0$, it is compact as well.
	Thus $g(e)\in\XHB$.  \sm

\item Assume that $e$ is such that $\zeta(e)=\infty$ but $a:=\liminf_{t\to\infty}e(t)<\infty$, and define
	$b:=\limsup_{t\to\infty} e(t)$. In the case $b>a$, $(T_e, r_e)$ is not Heine-Borel (and therefore not
	locally compact). Indeed, there is an $\eps>0$ with $a+3\eps<b$, and an increasing sequence $(t_n)$ in
	$\R_+$ with $e(t_n)\in [a+2\eps,a+3\eps]$ and $\inf_{u\in[t_n,t_{n+1}]} e(u) \le a+\eps$ for all
	$n\in\N$.  This means that $x_n:=\pi_e(t_n)$ defines a sequence of points in $\Bcl(\rho, a+3\eps)$ with
	mutual distances at least $\eps$.
	In the case $b=a$,
	\begin{equation}
		\mu_e\(\Bcl(\rho,b+1)\)=\lambda\bset{s\in\R_+}{e(s)\le b+1}=\infty,
	\end{equation}
	which means that $\mu_e$ is not boundedly finite. In both cases $g(e) \not\in \XHB$.
\end{enproof}\sm

Denote the space of continuous, transient excursions on $\R_+$ by
\begin{equation}
\label{e:expinfty}
   \E_{\mathrm{trans}}:=\bigl\{e\colon \R_+\to\R_+ \bigm| e \text{ is continuous},\, e(0)=0,\,
   \lim_{x\to\infty}e(x)=\infty\bigr\},
\end{equation}
and let for $e\in \E_{\mathrm{trans}}$ and $R>0$,
\begin{equation}
\label{e:lastexit}
   \xi_e(R):=\sup\{s\ge 0:\,e(s)<R\}<\infty
\end{equation}
denote the last visit to height $R>0$.

\begin{remark}[$\R$-trees under transient excursions]
	Let\/ $e\in\E_{\mathrm{trans}}$. Then  $g(e)$ is a Heine-Borel boundedly finite measure \Rtree\
	with precisely one end at infinity, i.e., there is a unique isometry $\vphi\colon[0,\infty)\to T_e$ with
	$\vphi(0)=\rho_e$. Indeed, the map $\vphi_e:=\pi_e\circ \xi_e$ is such an isometry.
	Assume that $\psi$ is a further such isometry and fix $R>0$. We show that $\psi(R)=\vphi_e(R)$.
	Choose $t\in \R_+$ with $\pi_e(t)=\psi(R)$. Because $\psi$ is an isometry, we have $e(t)=R$, and
	consequently $t\le \xi_e(R)$. Choose $S > \sup_{u\in[0,\xi_e(R)]}e(u)$ and $s\in\pi_e^{-1}\(\psi(S)\)$.
	Note that $s>\xi_e(R)$ and $e(s)=S$. Therefore $S-R = r_e\(\psi(S), \psi(R)\) = S+R-2\inf_{u\in[t,s]}
	e(u)$, and hence $\inf_{u\in [t,\xi_e(R)]} e(u) \ge \inf_{u\in[t,s]} e(u) = R$.
	This implies $r_e\(\psi(R), \vphi_e(R)\) = 2R-\inf_{u\in[t,\xi_e(R)]}e(u) = 0$.
\label{Rem:007}
\hfill$\qed$
\end{remark}\sm

\begin{proposition}[continuity of glue map]
	The glue map $g\colon \E_{\mathrm{trans}}\to\XHB$ is continuous if\/ $\E_{\mathrm{trans}}$ is equipped with the
	topology of uniform convergence on compacta, and\/ $\XHB$ with the Gromov-Hausdorff-vague topology.
\label{P:glue}
\end{proposition}\sm

\begin{proof} Let $(e_n)_{n\in\mathbb{N}}$ and $e$  in $\Etran$ be such that $e_n\tno e$ uniformly on compacta.
Put $\mathbb{W}_+:=\{R\ge 0:\,\lambda\{s\ge 0:\,e(s)=R\}>0\}$. Standard arguments show that $\mathbb{W}_+$ is at most countable.
Recall $\xi_e(R)$ from (\ref{e:lastexit}) and note that for all $R>0$, the map $e\mapsto \xi_e(R)$ is continuous with
respect to the uniform topology on compacta. Thus for all $R\in[0,\infty)\setminus\mathbb{W}_+$,
\begin{equation}
\label{e:gluecon}
   e_n\restricted{[0,\xi_{e_n}(R)]}\tno e\restricted{[0,\xi_{e}(R)]},
\end{equation}
which in turn implies that $g(e_n)\restricted{R}\to g(e)\restricted{R}$ Gromov-Hausdorff-weakly (see, for example, \cite[Proposition~2.9]{AbrahamDelmasHoscheit:exittimes}). Therefore $g(e_n)\tno g(e)$ Gromov-Hausdorff-vaguely by Definition~\ref{Def:002}.
\end{proof}\sm

We illustrate the usefulness of \propref{glue} with an example about the scaling limit of a size-biased branching tree (compare \cite{Kallenberg1977,GorostizaWakolbinger1991}
for a probabilitistic representation of this tree).
\begin{example}[Kallenberg-Kesten tree]
	Consider a (discrete time) Galton-Watson tree with a finite variance, mean~$1$ offspring distribution
	$p=(p_n)_{n\in\mathbb{N}_0}$.
	Let\/ $T'$ be the so-called \emph{Kallenberg-Kesten tree}, which is a random graph theoretic tree that is distributed
	like this Galton-Watson tree conditioned on survival.
	The simple, nearest neighbor random walk on $T'$, and scaling limits thereof, are of interest because of the
	``subdiffusive'' behavior (see \cite{MR88b:60232, BarlowKumagai2006}).
	The random walk is associated to the degree measure, defined in \exref{002}, as ``speed measure'' (see
	\cite[Section~7.4]{AthreyaLoehrWinter}).
	As in \exref{002}, we construct the (equivalent) rooted, measured \Rtree\ $(T, r, \rho, \mudegT)$, corresponding to $T'$.
	In the particular case of a geometric offspring distribution, i.e., $p_n:=2^{-(1+n)}$ for all\/
	$n\in\mathbb{N}_0$, we can code the tree with the length measure instead of the degree measure as follows:
	\begin{equation}
		(T,r,\rho,\lambdaT) \eqlaw g(\tilde{W}),
	\end{equation}
	where $\eqlaw$ denotes equivalence in law and, for all $t\ge 0$,
	$\tilde{W}_t:=W_t-2\inf_{s\in[0,t]}W_s$,
	with a simple random walk path $(W_n)_{n\in\mathbb{N}}$ linearly interpolated.
	We refer to $\TKalgeo:=(T,r,\rho, \mudegT)$ as the \emph{discrete Kallenberg
	tree} with geometric offspring distribution.

	As $W$ converges, after Brownian rescaling, weakly in path space towards standard Brownian motion $B$, we have
	\begin{equation}
	\label{e:WtoB}
	   \big(n^{-1}\tilde{W}_{n^2t}\big)_{t\ge 0}\Tno\big(\tilde{B}_t\big)_{t\ge 0},
	\end{equation}
	where $\tilde{B}_t:=B_t-2\inf_{s\in[0,t]}B_s$. It is shown in \cite{Pitman1975} that $(\tilde{B}_t)_{t\ge 0}$
	equals in law the unique strong solution of the stochastic differential equation
	\begin{equation}
	\label{e:sde}
	   X_t:=\tfrac{1}{X_t}\mathrm{d}t+\mathrm{d}B_t,\;t> 0,   \qquad X_0=0.
	\end{equation}
	Note that this solution is a three dimensional Bessel process (i.e.\ the radial path of a three dimensional
	Brownian motion). We refer to $g(X)$ as the \emph{continuum Kallenberg-Kesten tree}, $\TKal$.

	Because, almost surely, a realization $e:= n^{-1}\Wt_{n^{2}\bcdot}$ has slope
	$\pm n$ almost everywhere, we have $\mu_e=n^{-1}\lambdaT[(T_e, r_e)]$.
	Hence, by  Proposition~\ref{P:glue}, \eqref{e:WtoB} implies
	\begin{equation}
	\label{e:013}
	   (T, n^{-1}r, \rho, n^{-2}\lambdaT)\eqlaw g\(n^{-1}\tilde{W}_{n^{2}\boldsymbol{\cdot}}\)
	   \,\Tno\,
	   g(X)=\TKal,
	\end{equation}
	Gromov-Hausdorff-vaguely. By Example~\ref{Exp:002}, this also implies
	\begin{equation} \label{e:012}
		(T, n^{-1}r, \rho, n^{-2}\mudegT)
			\Tno
		\TKal,
	\end{equation}
	Gromov-Hausdorff-vaguely. In words, if we consider the discrete Kallenberg-Kesten tree and rescale
	the edge length to become $n^{-1}$, and then equip it with the measure which assigns mass
	$\tfrac{1}{2}n^{-2}\deg(x)$ to each branch point $x$, then this discrete measure tree converges
	weakly with respect to the Gromov-Hausdorff-vague topology to the continuum Kallenberg-Kesten tree.
	This implies that the simple, nearest neighbor random walk on the rescaled discrete Kallenberg-Kesten tree converges,
	if sped up by a factor of $n^3$, to the Brownian motion on $\TKal$, according to Theorem~1 of\/ \cite{AthreyaLoehrWinter}.
	See also Section~7.4 there.
\label{Exp:003}
\hfill$\qed$
\end{example}\sm

\bibliographystyle{alpha}
\bibliography{lit}

\begin{thebibliography}{ABBGM13}

\bibitem[ABBGM13]{AddarioBroutinGoldschmidtMiermont}
Lougi Addario-Berry, Nicolas Broutin, Christina Goldschmidt, and Gr{\'e}gory
  Miermont.
\newblock The scaling limit of the minimum spanning tree of the complete graph.
\newblock arXiv:1301.1664, 2013.

\bibitem[ADH13]{AbrahamDelmasHoscheit2013}
Romain Abraham, Jean-Fran{\c{c}}ois Delmas, and Patrick Hoscheit.
\newblock A note on the {G}romov-{H}ausdorff-{P}rokhorov distance between
  (locally) compact metric measure spaces.
\newblock {\em Electron. J. Probab.}, 18(14):1--21, 2013.

\bibitem[ADH14]{AbrahamDelmasHoscheit:exittimes}
Romain Abraham, Jean-Fran{\c{c}}ois Delmas, and Patrick Hoscheit.
\newblock Exit times for an increasing {L\'e}vy tree-valued process.
\newblock {\em Probab. Theory Related Fields}, 159(1-2):357--403, 2014.

\bibitem[Ald93]{Aldous1993}
David Aldous.
\newblock The continuum random tree {III}.
\newblock {\em Ann. Probab.}, 21:248--289, 1993.

\bibitem[ALW15]{AthreyaLoehrWinter}
Siva Athreya, Wolfgang L{\"o}hr, and Anita Winter.
\newblock Invariance principle for variable speed random walks on trees.
\newblock {\em Ann. Probab.}, in press, 49 pages, 2015.
\newblock arXiv:1404.6290.

\bibitem[BBI01]{BurBurIva01}
Dmitri Burago, Yuri Burago, and Sergei Ivanov.
\newblock {\em A course in metric geometry}, volume~33 of {\em Graduate studies
  in mathematics}.
\newblock AMS, Boston, MA, 2001.

\bibitem[BGMP14]{BlumbergGalMandellPancia}
Andrew~J. Blumberg, Itamar Gal, Michael~A. Mandell, and Matthew Pancia.
\newblock Robust statistics, hypothesis testing, and confidence intervals for
  persistent homology on metric measure spaces.
\newblock {\em Found. Comput. Math.}, 14(4):745--789, 2014.

\bibitem[BK06]{BarlowKumagai2006}
Martin~T. Barlow and Takashi Kumagai.
\newblock Random walk on the incipient infinite cluster on trees.
\newblock {\em Illinois J. Math.}, 50(1):33--65, 2006.

\bibitem[Car14]{Carlson2014}
Gunnar Carlson.
\newblock Topological pattern recognition for point cloud data.
\newblock {\em Acta Numerics}, 23:289--368, 2014.

\bibitem[CH13]{CurienHaas2013}
Nicolas Curien and B{\'en\'e}dicte Haas.
\newblock The stable trees are nested.
\newblock {\em Probab. Theory Related Fields}, 157:847--883, 2013.

\bibitem[Coh80]{Cohn1980}
Donald~L. Cohn.
\newblock {\em Measure Theory}.
\newblock Birkh{\"a}user, 1980.

\bibitem[Cro08]{Croydon08}
David~A. Croydon.
\newblock Volume growth and heat kernel estimates for the continuum random
  tree.
\newblock {\em Probab. Theory Related Fields}, 140:207--238, 2008.

\bibitem[DLG02]{MR1954248}
Thomas Duquesne and Jean-Fran{\c{c}}ois Le~Gall.
\newblock Random trees, {L\'e}vy processes and spatial branching processes.
\newblock {\em Ast\'erisque}, 281:vi+147, 2002.

\bibitem[Duq03]{Duquesne2003}
Thomas Duquesne.
\newblock A limit theorem for the contour process of conditioned
  {G}alton-{W}atson trees.
\newblock {\em Ann. Probab.}, 31:996--1027, 2003.

\bibitem[DW14]{DuquesneWang14}
Thomas Duquesne and Guanying Wang.
\newblock Exceptionally small balls in stable trees.
\newblock {\em Bull. Soc. Math. France}, 142(2):223--254, 2014.

\bibitem[EW06]{EvansWinter2006}
Steven~N. Evans and Anita Winter.
\newblock Subtree prune and re-graft: A reversible real-tree valued {M}arkov
  chain.
\newblock {\em Ann. Probab.}, 34(3):918--961, 2006.

\bibitem[FOT11]{FukushimaOshimaTakeda1994}
Masatoshi Fukushima, Yoichi Oshima, and Masayoshi Takeda.
\newblock {\em Dirichlet Forms and Symmetric {M}arkov Processes}, volume~19 of
  {\em de Gruyter Studies in Mathematics}.
\newblock Walter de Gruyter \& Co., Berlin, second edition, 2011.

\bibitem[Fuk87]{Fukaya1987}
Kenji Fukaya.
\newblock Collapsing of {R}iemannian manifolds and eigenvalues of {L}aplace
  operators.
\newblock {\em Invent. Math.}, 87:517--547, 1987.

\bibitem[GPW09a]{GrevenPfaffelhuberWinter2009}
Andreas Greven, Peter Pfaffelhuber, and Anita Winter.
\newblock Convergence in distribution of random metric measure spaces
  ({$\Lambda$}-coalescent measure trees).
\newblock {\em Probab. Theory Related Fields}, 145(1-2):285--322, 2009.

\bibitem[GPW09b]{GrevenPopovicWinter2009}
Andreas Greven, Lea Popovic, and Anita Winter.
\newblock Genealogy of catalytic branching models.
\newblock {\em Annals of Applied Probability}, 19(3):1232--1272, 2009.

\bibitem[GPW13]{GrevenPfaffelhuberWinter2013}
Andreas Greven, Peter Pfaffelhuber, and Anita Winter.
\newblock Tree-valued resampling dynamics: {M}artingale {P}roblems and
  applications.
\newblock {\em Probab. Theory Related Fields}, 155(3--4):789--838, 2013.

\bibitem[Gro99]{MR2000d:53065}
Misha Gromov.
\newblock {\em Metric structures for {R}iemannian and non-{R}iemannian spaces},
  volume 152 of {\em Progress in Mathematics}.
\newblock Birkh\"auser Boston Inc., Boston, MA, 1999.
\newblock Based on the 1981 French original.

\bibitem[GW91]{GorostizaWakolbinger1991}
Luis~G. Gorostiza and Anton Wakolbinger.
\newblock Persistence criteria for a class of critical branching particle
  systems in continuous time.
\newblock {\em Ann. Probab.}, 19(1):266--288, 1991.

\bibitem[HM12]{HaasMiermont2012}
B{\'e}n{\'e}dicte Haas and Gr{\'e}gory Miermont.
\newblock Scaling limits of {M}arkov branching trees with applications to
  {G}alton-{W}atson and random unordered trees.
\newblock {\em Ann. Probab.}, 40(6):2589--2666, 2012.

\bibitem[HW14]{HeWinkel}
Hui He and Matthias Winkel.
\newblock Invariance principles for pruning processes of {Galton-Watson} trees.
\newblock arXiv:1409.1014, 2014.

\bibitem[Kal77]{Kallenberg1977}
Olav Kallenberg.
\newblock Stability of critical cluster fields.
\newblock {\em Math. Nachr.}, 77:7--43, 1977.

\bibitem[Kes86]{MR88b:60232}
Harry Kesten.
\newblock Subdiffusive behavior of random walk on a random cluster.
\newblock {\em Ann. Inst. H. Poincar\'e Probab. Statist.}, 22(4):425--487,
  1986.

\bibitem[KS03]{KuwaeShioya2003}
Kazhiro Kuwae and Takashi Shioya.
\newblock Convergence of spectral structure: a functional analytic theory and
  its application to spectral geometry.
\newblock {\em Analysis and Geometry}, 11(4):599--673, 2003.

\bibitem[LC57]{LeCam1957}
Lucien Le~Cam.
\newblock Convergence in distribution of stochastic processes.
\newblock {\em Univ. California Publ. Statist.}, 2:207--236, 1957.

\bibitem[LG07]{LeGall2007}
Jean-Fran{\c{c}}ois Le~Gall.
\newblock The topological structure of scaling limits of large planar maps.
\newblock {\em Invent. Math.}, 169(3):621--670, 2007.

\bibitem[L{\"o}h13]{Loehr}
Wolfgang L{\"o}hr.
\newblock Equivalence of {G}romov-{P}rohorov- and {G}romov's
  {$\underline\square_\lambda$}-metric on the space of metric measure spaces.
\newblock {\em Electron. Commun. Probab.}, 18(17):1--10, 2013.

\bibitem[LVW15]{LoehrVoisinWinter2013}
Wolfgang L{\"o}hr, Guillaume Voisin, and Anita Winter.
\newblock Convergence of bi-measure {$\mathbb R$-trees} and the pruning
  process.
\newblock {\em Ann. Inst. H. Poincar\'e Probab. Statist.}, 51(4):1342--1368,
  2015.

\bibitem[Mie09]{Miermont2009}
Gr{\'e}gory Miermont.
\newblock Tessellations of random maps of arbitrary genus.
\newblock {\em Ann. Sci. {\'E}c. Norm. Sup{\'e}r. (4)}, 42(5):725--781, 2009.

\bibitem[MM11]{MarckertMiermont2011}
Jean-Fran{\c{c}}ois Marckert and Gr{\'e}gory Miermont.
\newblock The {CRT} is the scaling limit of unordered binary trees.
\newblock {\em Random Structures Algorithms}, 38(4):467--501, 2011.

\bibitem[Pit75]{Pitman1975}
James~W. Pitman.
\newblock {One-dimensional Brownian motion and the three-dimensional Bessel
  process}.
\newblock {\em Adv. in Appl. Probab.}, 7:511--526, 1975.

\bibitem[Vil08]{Villani2008}
C\'{e}dric Villani.
\newblock Optimal transport, old and new.
\newblock In {\em {\'Ecole d'\'Et\'e} de {Probabilit\'es} de Saint Flour
  XXXV-2005}, volume 1920 of {\em Lecture Notes in Mathematics}.
  Springer-Verlag, 2008.

\end{thebibliography}

\end{document}